\documentclass[12pt,oneside,reqno]{amsart}
\usepackage{amsmath,amssymb,amsthm}
\usepackage[margin=1.1in]{geometry}
\usepackage{booktabs}
\usepackage[hidelinks]{hyperref}
\hypersetup{pdftitle={Edge complexity of complete bipartite and complete multipartite graphs},pdfauthor={S. Ferguson, V. Gupta, A. Iosevich, J. Iosevich, C. Kelly, G. Smith, S. Yang}}
\newtheorem{theorem}{Theorem}[section]
\newtheorem{lemma}[theorem]{Lemma}
\newtheorem{proposition}[theorem]{Proposition}
\newtheorem{corollary}[theorem]{Corollary}
\theoremstyle{definition}
\newtheorem{definition}[theorem]{Definition}
\newtheorem{question}[theorem]{Question}
\theoremstyle{remark}
\newtheorem{remark}[theorem]{Remark}

\title[Edge complexity of complete multipartite graphs]{Edge complexity of complete bipartite and complete multipartite graphs}
\author{S. Ferguson}
\address{Department of Mathematics, New York University}
\email{sjf370@nyu.edu}

\author{V. Gupta}
\address{Department of Mathematics, University of Rochester}
\email{vishalgupta@rochester.edu}

\author{A. Iosevich}
\address{Department of Mathematics, University of Rochester}
\email{alex.iosevich@rochester.edu}

\author{J. Iosevich}
\address{Department of Applied Mathematics,
Rochester Institute of Technology}
\email{joshuaiosevich@gmail.com}

\author{C. Kelly}

\author{G. Smith}

\author{S. Yang}
\thanks{The third listed author was supported in part by NSF grant DMS-2506858.}
\date{October 26, 2026}

\begin{document}
\begin{abstract}
We study edge complexity, the minimum Fourier $\ell^1/\ell^2$
ratio of a graph's adjacency matrix over cyclic vertex labelings.
For $K_{m,n}$, with $1\leq m\leq n$, we reduce the problem to
minimizing the Fourier $\ell^1$-norm of indicators of $m$-element
subsets of $\mathbb Z_{m+n}$. We obtain the sharp lower bound
\[
\operatorname{FR}_{\min}(K_{m,n})
\geq\frac{4n^2}{(m+n)\sqrt{2mn}},
\]
with equality precisely when $m\mid n$, and classify all
labelings attaining it. We determine all possible Fourier ratios and minimizing
labelings for $K_{2,n}$, and all minimizing labelings for
$K_{3,n}$. For the latter, put $N=n+3$. The optimal
three-element sets are precisely the cosets of the subgroup
of order $3$ in $\mathbb Z_N$ when $3\mid N$, and precisely
the arithmetic progressions with invertible common difference
otherwise. The three-element classification uses polynomial
minorants and additive energies.

We prove that attainment of the graph-energy bound is invariant
under uniform independent set blow-ups. Edge complexity itself
can decrease by an arbitrarily large factor under a twofold
blow-up, and we also give an explicit finite example.
For complete multipartite graphs, we obtain an exact labeling
formula and a sharp lower bound, and classify all part sizes
and labelings attaining that bound. This bound strictly improves
the graph-energy bound whenever the part sizes are unequal.
We also prove a lower bound for arbitrary simple graphs in
terms of their average degree and maximum absolute degree
deviation.

For Laplacian projectors, we apply the equality criterion for the multiplicity bound to unbalanced
complete bipartite graphs, whose top-projector minimum reduces
to the same set minimization. Finally, we apply existing Fourier
recovery and sampling theorems to approximate a graph from random
adjacency samples under a fixed known labeling with bounded
Fourier ratio, obtaining a bound on the number of incorrectly
recovered edges.
\end{abstract}
\maketitle
\tableofcontents

\section{Introduction}

The adjacency matrix of a simple undirected graph records which pairs of vertices are joined by an
edge. Although its eigenvalues do not depend on the order in which the vertices are listed, its Fourier transform does depend on this order. We are
interested in how to label a graph's vertices so that its adjacency matrix has
the smallest possible Fourier ratio, a measure of Fourier complexity and compressibility. In this paper, we study this
question for complete bipartite and complete multipartite graphs.
We determine all minimizing labelings when one bipartite part has
two or three vertices, and obtain sharp lower bounds with complete
equality classifications for arbitrary part sizes.
We also investigate
how the minimum changes under blow-ups, how the same labeling problem
arises for Laplacian projectors, and how a small Fourier ratio allows
us to recover edge information from partial observations of a graph.

Bipartite networks encode relations between two types of
objects, such as authors and papers or stores and warehouses.
Larremore, Clauset, and Jacobs
\cite{LarremoreClausetJacobs2014} show that projecting a
bipartite network onto one vertex class can lose information
relevant to community structure, and develop a stochastic
block model that retains both vertex types. This provides
one motivation for recovering the original incidence
information when some adjacencies are unobserved.

Complete bipartite and complete multipartite graphs provide
explicit benchmarks for studying how the arrangement and
sizes of vertex classes affect Fourier complexity. When the
parts and completeness are known, the graph itself is already
determined. The recovery question concerns incomplete
networks; Section~\ref{sec:recovery} treats this question
under a known labeling with a bound for the Fourier ratio
of the actual adjacency function.

\subsection{Edge complexity and the bipartite reduction}

Let $G=(V,E)$ be a finite simple undirected graph with $N$ vertices
and $s>0$ edges. A labeling is a bijection $\pi:\mathbb Z_N\to V$.
It identifies the adjacency matrix with the function
\[
f_\pi(x,y)=1_{\{\pi(x)\sim\pi(y)\}},
\qquad (x,y)\in\mathbb Z_N^2.
\]
For a complex-valued function $f$ on $\mathbb Z_N^d$, we use the unitary Fourier
transform
\[
\widehat f(\xi)=\frac{1}{N^{d/2}}
\sum_{x\in\mathbb Z_N^d}f(x)e^{-2\pi i x\cdot\xi/N}.
\]
All norms use counting measure. In particular, matrix norms are
entrywise norms: the $\ell^2$-norm is the Frobenius norm.
Parseval's identity reads $\|\widehat f\|_2=\|f\|_2$, and for
$f\neq0$ the Fourier ratio is
\[
\operatorname{FR}(f)=\frac{\|\widehat f\|_1}{\|\widehat f\|_2}.
\]
The Fourier ratio and its applications are studied in
\cite{AldalehEtAl2025,BursteinIosevichNathan2026}. Following Gupta
and Iosevich \cite{GuptaIosevich2026}, we define the edge complexity by
\[
\operatorname{FR}_{\min}(G)
=\min_{\pi:\mathbb Z_N\to V}\operatorname{FR}(f_\pi).
\]
There are finitely many labelings, so this minimum is attained.
The assumption $s>0$ excludes the zero adjacency matrix, for which
the ratio is undefined. Since each undirected edge contributes two
entries equal to one, $\|f_\pi\|_2=\sqrt{2s}$ for every labeling.

There is a general lower bound in terms of the adjacency spectrum.
If $\lambda_1(G),\ldots,\lambda_N(G)$ are the adjacency eigenvalues,
the graph energy is
\[
\mathcal E(G)=\sum_{j=1}^N|\lambda_j(G)|.
\]
Gupta and Iosevich \cite{GuptaIosevich2026} proved that
\begin{equation}\label{eq:energy-bound}
\operatorname{FR}_{\min}(G)\geq\frac{\mathcal E(G)}{\sqrt{2s}}.
\end{equation}
Cyclic Cayley graphs attain this bound under their natural labelings.
We recall the equality criterion of \cite{GIIST2026}, with a proof,
in Section~\ref{sec:blowups}. Our results will also give families for
which a strictly stronger lower bound is available.

Related work considers both the computational and additive aspects
of edge complexity. Gupta and Iosevich
\cite{GuptaIosevichWeighted2026} proved inverse-polynomial
additive approximation hardness for weighted edge complexity
and established a periodic blow-up inequality.
Gupta et al.~\cite{GuptaEtAlAdditive2026} studied additive
edge energy, sparse Fourier approximation, and spectral
projector kernels.

For a complete bipartite graph, the labeling problem has a
particularly concrete form. Suppose $1\leq m\leq n$ and put
$N=m+n$. Assign an $m$-element set $C\subset\mathbb Z_N$ to one
part and its complement to the other. Rearranging the labels within
either part does not change the adjacency function, which we denote
by $f_C$. Put
\begin{equation}\label{eq:S-D-definitions}
S_C(t)=\sum_{c\in C}e^{-2\pi itc/N},
\qquad D_N(C)=\sum_{t=1}^{N-1}|S_C(t)|,
\end{equation}
and write
\[
d_{m,N}=\min_{\substack{C\subset\mathbb Z_N\\|C|=m}}D_N(C).
\]
The sum $S_C$ is unnormalized; in our convention it equals
$\sqrt N\widehat{1_C}$.

Our first result is the identity
\begin{equation}\label{eq:intro-formula}
\operatorname{FR}(f_C)=
\frac{2mn+2(n-m)D_N(C)+2D_N(C)^2}{N\sqrt{2mn}}.
\end{equation}
The right side is a strictly increasing function of $D_N(C)\geq0$.
Thus the two-dimensional graph problem reduces exactly to a
one-dimensional minimization over subsets of prescribed size.
Fourier inversion gives $D_N(C)\geq n$, with equality precisely
when $C$ is a coset of the subgroup of order $m$. Consequently,
\[
\operatorname{FR}_{\min}(K_{m,n})
\geq\frac{4n^2}{(m+n)\sqrt{2mn}},
\]
with equality if and only if $m\mid n$. We classify all labelings
attaining this bound. The equivalence between cosets and indicators
of Wiener norm one is classical; see the discussion preceding
Theorem~\ref{thm:coset}. Here it determines the equality cases in
the graph problem. When $m\nmid n$, no subgroup of order $m$
exists in $\mathbb Z_N$, and we must ask which other sets minimize
$D_N$.

\subsection{Littlewood's problem and additive structure}

This minimization problem has a classical ancestor. Let
$A\subset\mathbb Z$ consist of $m$ distinct integers, and put
\[
F_A(\theta)=\sum_{a\in A}e^{2\pi ia\theta},
\qquad \theta\in\mathbb R/\mathbb Z.
\]
We ask how small its average absolute value
\[
I(A)=\int_0^1|F_A(\theta)|\,d\theta
\]
can be. At zero all $m$ terms agree, but away from zero they can
cancel. Orthogonality fixes the $L^2$-norm at $\sqrt m$, independently
of the positions of the frequencies. The $L^1$-norm, on the other
hand, depends on how those frequencies are arranged.

Consecutive frequencies provide the guiding example. For
$A=\{0,\ldots,m-1\}$, the geometric-series formula gives
\[
|F_A(\theta)|=
\left|\frac{\sin(\pi m\theta)}{\sin(\pi\theta)}\right|,
\qquad \theta\notin\mathbb Z.
\]
This is the calculation underlying the Dirichlet kernel. Write
\[
J_m=2\int_0^{1/2}
\frac{|\sin(\pi m\theta)|}{\sin(\pi\theta)}\,d\theta
\]
for its integral over the circle. The initial interval $[0,1/m]$
contributes $O(1)$, since $|F_A|\leq m$. Farther from zero, the
numerator oscillates on a scale over which the denominator changes
comparatively little. Replacing the numerator by its period average
\[
m\int_0^{1/m}|\sin(\pi m\theta)|\,d\theta
=\int_0^1|\sin(\pi u)|\,du=\frac2\pi
\]
therefore suggests the approximation
\[
J_m\approx
2\int_{1/m}^{1/2}\frac{2/\pi}{\sin(\pi\theta)}\,d\theta
=\frac4{\pi^2}\log\cot\left(\frac{\pi}{2m}\right)
=\frac4{\pi^2}\log m+O(1).
\]
There is also a geometric way to see the logarithm. Near zero,
this averaging models the integrand by $2/(\pi^2\theta)$.
Each interval $[r,2r]$ then contributes approximately
$(2/\pi^2)\log2$. There are $\log_2m+O(1)$ such intervals between
$1/m$ and $1/2$, and symmetry accounts for the other half of the
circle. These are heuristic explanations. Averaging over successive
periods gives the rigorous classical asymptotic
\[
J_m=\frac4{\pi^2}\log m+O(1);
\]
see \cite[Section~1.2]{BloomGreen2026}.

The logarithmic lower bound has a short direct proof. For $m\geq2$,
put $q=\lfloor m/2\rfloor$. Using $\sin(\pi\theta)\leq\pi\theta$
and substituting $u=m\theta$, we obtain
\[
\begin{aligned}
J_m
&\geq\frac2\pi\int_0^{m/2}\frac{|\sin(\pi u)|}{u}\,du\\
&\geq\frac2\pi\sum_{j=1}^q\frac1j
       \int_{j-1}^j|\sin(\pi u)|\,du
=\frac4{\pi^2}\sum_{j=1}^q\frac1j.
\end{aligned}
\]
Here we discard the possible final partial interval, use
$1/u\geq1/j$ on $(j-1,j]$, and integrate one period of the absolute
sine. Since
\[
\sum_{j=1}^q\frac1j\geq\log(q+1),
\qquad q+1\geq\frac{m+1}{2}\geq\sqrt m,
\]
we have $J_m\geq(2/\pi^2)\log m$.

Littlewood asked whether every set of $m$ distinct integer
frequencies must have an $L^1$-norm bounded below by a positive
absolute constant times $\log m$. The question appears in Hardy
and Littlewood's 1948 paper on rearrangements
\cite{HardyLittlewood1948}. The consecutive-frequency example
supplies both the proposed lower bound and the reason its order
of growth cannot be improved. The conjecture was proved
independently in 1981 by McGehee, Pigno, and Smith
\cite{McGeheePignoSmith1981} and by Konyagin \cite{Konyagin1982}:
there is an absolute constant $c>0$ such that
\[
I(A)\geq c\log m
\]
for every set of $m\geq2$ distinct integers. Konyagin's paper
appeared in Russian in 1981; its cited English translation
appeared in 1982. This theorem determines the smallest order of
growth, but not the exact minimizing sets. The stronger question
of whether an arithmetic progression always gives the exact
minimum is usually called the strong Littlewood conjecture;
see \cite[Introduction]{GreenKonyagin2009}.

One reason to consider progressions comes from comparing a set
with a translate. For a positive integer $d$,
\[
(1-e^{2\pi id\theta})F_A(\theta)
=F_A(\theta)-F_{A+d}(\theta).
\]
Every common frequency cancels. If $A$ is a progression with
step $d$, only two endpoints remain. Conversely, each separate
chain of points spaced by $d$ contributes two endpoints, so a
single progression has the fewest possible surviving terms.
One might expect additional unmatched terms to increase
fluctuations and hence the average absolute value. This is only
a heuristic: recovering $F_A$ requires division by
$1-e^{2\pi id\theta}$, and its correlations with the remaining
terms prevent a direct comparison of the integrals.

The same calculation leads to additive energy. For $t\in\mathbb Z$,
let
\[
r_A(t)=|A\cap(A+t)|
=|\{(a,b)\in A^2:a-b=t\}|.
\]
There are $2(m-r_A(d))$ surviving terms in $F_A-F_{A+d}$.
Orthogonality therefore gives the exact identity
\[
\int_0^1|(1-e^{2\pi id\theta})F_A(\theta)|^2\,d\theta
=2(m-r_A(d)).
\]
A large overlap with a translate means a small mean-square
remainder after cancellation. Additive energy collects these
overlaps over all translates:
\[
\begin{aligned}
E_2(A)
&=|\{(a,b,c,e)\in A^4:a+b=c+e\}|\\
&=\sum_{t\in\mathbb Z}r_A(t)^2
=\int_0^1|F_A(\theta)|^4\,d\theta.
\end{aligned}
\]
The counting identity follows by rewriting $a+b=c+e$ as
$a-c=e-b$; the integral identity follows by expansion and
orthogonality. For a progression with step $d$,
$r_A(kd)=m-|k|$ when $|k|<m$, with all other overlaps zero. Hence
\[
E_2(A)=m^2+2\sum_{j=1}^{m-1}j^2=\frac{2m^3+m}{3}.
\]

There is a rigorous necessary condition connecting this energy
to a small $L^1$-norm. H\"older's inequality with exponents $3/2$
and $3$ gives
\[
m=\int_0^1|F_A(\theta)|^{2/3}|F_A(\theta)|^{4/3}\,d\theta
\leq I(A)^{2/3}E_2(A)^{1/3}.
\]
Thus
\[
E_2(A)\geq\frac{m^3}{I(A)^2},
\qquad
I(A)\leq K\log m
\ \Longrightarrow\
E_2(A)\geq\frac{m^3}{K^2(\log m)^2};
\]
see \cite[Introduction]{BloomGreen2026}.
Small average amplitude requires many additive relations.
Nevertheless, this does not say that increasing additive energy
always decreases the $L^1$-norm. For a progression, the equivalent
bound $I(A)\geq m^{3/2}/E_2(A)^{1/2}$ gives only a constant-order
lower bound, not the logarithmic growth found above.

Our graph problem replaces integration over the circle by averaging
over its $N$ equally spaced points. The relevant Wiener norm is
\[
\|1_C\|_{A(\mathbb Z_N)}
=\frac1N\sum_{t\in\mathbb Z_N}|S_C(t)|
=\frac{m+D_N(C)}N.
\]
Minimizing this norm for fixed $m,N$ is therefore equivalent to
minimizing $D_N(C)$. The change of group matters. A subgroup coset
has Wiener norm one regardless of its size, so no logarithmic
lower bound in $|C|$ can hold uniformly over all finite cyclic
groups. At prime orders, the absence of nontrivial proper subgroups
removes this obstruction. Green and Konyagin
\cite{GreenKonyagin2009} studied the finite cyclic problem, and
Sanders \cite{Sanders2007} obtained stronger lower bounds. We use
these estimates in our blow-up examples. Our small-part
classifications address the complementary question of exact
minimizers.

\subsection{Two- and three-element partite sets}

For two labels, write $C=\{0,d\}$ after translation and put
$g=\gcd(d,N)$ and $q=N/g$. The Fourier sum depends only on the
additive order $q$ of the difference. We evaluate the sums in
terms of cotangents and cosecants, compare the even and odd
branches, and determine all possible Fourier ratios of $K_{2,n}$
for $n\geq2$. Exactly $\tau(N)-1$ values occur, where $\tau(N)$
is the number of positive divisors of $N$. If $N$ is even, the
minimizers are precisely the antipodal pairs; if $N$ is odd,
they are precisely the pairs whose difference is invertible
modulo $N$. The trigonometric evaluations are classical
Wiener-norm formulas; Section~\ref{sec:two} recalls their
provenance and gives a direct derivation and the comparisons
needed for the graph problem.

For three labels, we obtain a complete classification at every
order. If $N\geq4$, the minimizing three-element subsets of
$\mathbb Z_N$ are precisely
\[
\begin{cases}
\{t,t+N/3,t+2N/3\},&3\mid N,\\
\{t,t+u,t+2u\},\quad u\in\mathbb Z_N^\times,&3\nmid N,
\end{cases}
\qquad t\in\mathbb Z_N.
\]
This follows from Theorems~\ref{thm:coset} and
\ref{thm:three-minimizers-every-order}, and determines all
minimizing labelings of $K_{3,n}$ for $n\geq1$, with $N=n+3$.
The invertibility condition is essential at composite orders:
a progression with three distinct elements need not have an
invertible common difference.

The proof uses additive energies to compare Fourier norms.
An even polynomial lying below the absolute value converts
moment identities into a lower bound for $D_N(C)$. The
three-term argument of mathworker21 and Fedja on MathOverflow
\cite{Mathworker21Fedja2022} supplies the starting method on the
circle. In the finite cyclic setting, extra additive relations
and common differences of different orders require further
comparisons. We use rational minorants of degrees six and eight.
The additional moment estimate is that $E_3(C)=93$ forces
$E_4(C)\leq667$ when $3\nmid N$. Indeed, the absence of
coincidences between distinct unordered three-term sums severely
restricts coincidences between four-term sums. When $3\nmid N$,
the resulting bounds exclude every nonprogression, while a separate
comparison by the additive order of the common difference identifies
the minimizing progressions. Here $E_j(C)$ counts equal sums of two
ordered $j$-tuples from $C$, with equality taken modulo $N$.

These classifications do not extend unchanged to larger parts.
In Subsection~\ref{sec:labelings}, we compare interval, balanced,
and arithmetic progression labelings. For $K_{7,9}$, an explicit
labeling improves on both the interval and balanced labelings,
even though $7$ and $16$ are relatively prime. Thus the
three-element result leaves a different problem beginning with
four-element sets, rather than settling the general minimization.

\subsection{Blow-ups and the role of the ambient group}

A uniform independent set blow-up $G^{(k)}$ replaces each vertex
by $k$ vertices with the same neighbors, and each original edge
by all edges between the corresponding two sets. Repeating a
labeling periodically preserves its Fourier ratio and gives
\[
\operatorname{FR}_{\min}(G^{(k)})
\leq\operatorname{FR}_{\min}(G).
\]
We prove that $G^{(k)}$ attains its graph-energy bound if and
only if $G$ does. The forward implication requires descending
an extremizing labeling of the enlarged graph; it does not
follow merely from the existence of a periodic lift.

Nevertheless, the edge complexity itself need not be preserved.
For an odd prime $p$, put
\[
G_p=K_{(p-1)/2,(p+1)/2}.
\]
For all sufficiently large $p$, we prove
\[
\operatorname{FR}_{\min}(G_p)
\geq c\frac{\log p}{(\log\log p)^3},
\]
where $c>0$ is absolute, whereas for every odd prime $p$,
\[
\operatorname{FR}_{\min}(G_p^{(2)})<5\sqrt2.
\]
Thus a twofold blow-up can decrease edge complexity by an
arbitrarily large factor. The lower bound uses Sanders's
improvement of the theorem of Green and Konyagin
\cite{GreenKonyagin2009,Sanders2007}. For the upper bound,
the new group $\mathbb Z_{2p}$ has a subgroup of index two;
removing one point from it gives a part of the required size
with bounded normalized Fourier norm. This explains why a
labeling of the enlarged graph can improve on every periodic lift.

We also give a finite example, certified by rational arithmetic:
\[
\operatorname{FR}_{\min}(K_{22,24})
<\operatorname{FR}_{\min}(K_{11,12}).
\]
Consequently, dividing both part sizes by a common divisor gives
an upper bound through periodic lifting, but does not in general
determine the minimum. These results motivate the questions about
repeated blow-ups and descent through cyclic quotients in
Section~\ref{sec:questions}.

\subsection{Complete multipartite graphs and a degree bound}

For more than two parts, a labeling is a partition
$\mathbb Z_N=C_1\cup\cdots\cup C_r$ with $|C_j|=n_j$.
The adjacency function is
\[
f(x,y)=1-\sum_{j=1}^r1_{C_j}(x)1_{C_j}(y).
\]
Its Fourier transform gives an exact labeling formula
(Theorem~\ref{thm:multipartite-formula}). Unlike the bipartite
case, the Fourier sums of the parts do not reduce to a single
sum and its negative. Their phases affect the answer, so
minimizing the separate indicator norms does not by itself
solve the partition problem.

Despite this interaction, we obtain a sharp bound with a complete
equality classification. Let
\[
N=\sum_{j=1}^r n_j,\qquad
\mu=\frac1N\sum_{j=1}^r n_j^2,\qquad
a=\min_j n_j,\qquad 2s=N(N-\mu).
\]
Then
\[
\operatorname{FR}_{\min}(K_{n_1,\ldots,n_r})
\geq\frac{2(N+\mu-2a)}{\sqrt{2s}}.
\]
Equality holds precisely for balanced part sizes, or, after
reordering, for $h=r-1$ parts of size $a$ and one part of size
$kha$, where $k\geq1$ is an integer. In the balanced case the
extremizing parts are the cosets of the subgroup of order $a$.
In the unbalanced case, take the subgroups $H\subseteq L$ of
orders $a$ and $ha$. The small parts are exactly the $h$ cosets
of $H$ in one translate of $L$, and the large part is its
complement. These are all extremizing labelings
(Theorem~\ref{thm:multipartite-sharp-equality}).

For example, the tripartite graphs $K_{a,a,qa}$ attain the bound
for every positive even integer $q$, giving
\[
\operatorname{FR}_{\min}(K_{a,a,qa})
=\frac{4(q^2+q+1)}{(q+2)\sqrt{4q+2}}.
\]
The multipartite lower bound strictly improves the graph-energy
bound whenever the part sizes are unequal. Its equality families
also preserve edge complexity under every uniform blow-up,
including unbalanced examples that do not attain the energy
bound. Outside the equality families, the exact multipartite
minimum remains undetermined in general.

A related argument gives a lower bound for every simple graph
with $s>0$ edges. If $\overline d=2s/N$ and
$\delta=\max_x|\deg(x)-\overline d|$, then
\[
\operatorname{FR}_{\min}(G)
\geq\frac{2\overline d+4\delta}{\sqrt{2s}}.
\]
This estimate uses the zero diagonal of the adjacency matrix
and Fourier inversion of the degree function. For a complete
multipartite graph, put $b=\max_j n_j$. The estimate strengthens
the preceding bound when $b-\mu>\mu-a$ and coincides with it
otherwise. Its equality problem is separate from the
multipartite equality classification above.

\subsection{Laplacian projectors}

The labeling problem also applies to matrices other than
adjacency matrices. For the combinatorial Laplacian
$L=\operatorname{diag}(\deg(v))-A$, let $\Pi_\lambda$ be the
orthogonal projector onto the $\lambda$-eigenspace. We regard
its kernel as a function of two vertex labels and minimize its
Fourier ratio in the same way. Gupta and Iosevich
\cite{GuptaIosevich2026} proved the multiplicity bound
\[
\operatorname{FR}(\Pi_\lambda)
\geq\sqrt{\operatorname{rank}\Pi_\lambda}
\]
and showed that a circulant graph attains all these bounds
under its natural labeling. We recall the equality criterion of
\cite[Proposition~11.2]{GuptaEtAlAdditive2026}:
under a fixed cyclic labeling, a nonzero orthogonal projector
attains the multiplicity bound if and only if it is circulant.

For $K_{m,n}$ with $1\leq m<n$ and $N=m+n$, the constant
projector is the only nonzero spectral projector attaining its
multiplicity bound. The projector onto the top Laplacian
eigenspace satisfies the more precise identity
\[
\operatorname{FR}_{\min}(\Pi_N)
=\frac{d_{m,N}^2}{mn}\geq\frac nm>1,
\]
with equality in the lower bound $n/m$ exactly when $m\mid n$.
Its minimizing labelings are exactly those minimizing adjacency.
Thus the bipartite set problem determines both minima through
different increasing functions of $D_N(C)$.

For comparison, we recall the result of Gupta and Iosevich that
the parity labeling of $K_{m,m}$ simultaneously attains the
multiplicity bound for every nonzero Laplacian projector.
In the unbalanced case, the exact intermediate-projector minima
remain open in general. Their formulas retain Fourier phases
as well as magnitudes, and their minimizing labelings need not
minimize adjacency. Section~\ref{sec:questions} gives an example
and asks when common minimizing labelings exist.

\subsection{Approximate recovery from adjacency samples}

One motivation for finding labelings with small Fourier ratio
is recovery from partial observations. We consider a bipartite
graph whose vertices represent stores and warehouses, with an
edge when a delivery can be made within one day. Such a graph
need not be complete. Its Fourier ratio therefore requires its
own bound; the complete bipartite calculation does not supply
one merely from the part sizes or their divisibility.

Suppose a known labeling is fixed before sampling and its
adjacency function satisfies $\operatorname{FR}(f)\leq R$.
Choose $q$ distinct ordered pairs uniformly from
$\mathbb Z_N^2$, observe their adjacency values, and minimize
$\|\widehat g\|_1$ subject to agreement at those observations.
This is Fourier basis pursuit \cite{ChenDonohoSaunders1998}.
We apply the recovery principle of Burstein, Iosevich, and Nathan
\cite{BursteinIosevichNathan2026}, together with the improved
sampling theorem of Burstein, Iosevich, and Krause
\cite{BursteinIosevichKrause2026}.

Theorem~\ref{thm:graph-approximation} bounds the relative
$\ell^2$ error. Symmetrizing the approximation and rounding
at $1/2$ then converts it into a graph. If
$\|g-f\|_2\leq\varepsilon\sqrt{2s}$, the recovered edge set
$E_*$ satisfies
\[
|E\triangle E_*|\leq4\varepsilon^2s,
\]
counting both inserted and omitted edges. More explicitly, for
$0<\delta<1$ put $k_\delta=\lceil16R^2/\delta\rceil$. For an
absolute constant $C>0$, if $k_\delta\leq N^2$ and
\[
Ck_\delta\log\left(\frac{2eN^2}{k_\delta}\right)
\log(2k_\delta)\leq q\leq\frac{N^2}{2},
\]
then $|E\triangle E_*|\leq\delta s$ with probability at least
$1-CN^{-2}$. The sampling is from the full square, including
known zero entries, rather than only from the potentially
present store--warehouse edges. If the stated sample range is
unavailable, one can instead inspect all $mn$ possible edges.

This result approximates the labeled edge set under an available
labeling with a controlled ratio. A small minimum over all
labelings does not itself provide such a labeling. Nor does a
bounded ratio alone guarantee exact recovery from few queries:
deleting one unknown edge from a complete bipartite graph gives
an obstruction. These distinctions lead to further questions
about obtaining useful labelings from geometry or partial data,
exploiting the known bipartition when sampling, and determining
the sharp dependence on the ratio and the allowed number of errors.

The paper is organized as follows. Section~\ref{sec:bipartite}
proves the bipartite formula and its sharp lower bound,
then compares interval, balanced, and arithmetic progression
labelings. Section~\ref{sec:two} determines all Fourier
ratios and minimizing labelings for $K_{2,n}$.
Section~\ref{sec:three} gives the complete classification
of three-element minimizing sets.
Section~\ref{sec:blowups} studies blow-ups and the
graph-energy bound, including a certified finite example of
strict decrease in Subsection~\ref{sec:certificate}.
Section~\ref{sec:multipartite} proves the multipartite results
and the degree bound, and Section~\ref{sec:harmonic} treats
Laplacian projectors. Section~\ref{sec:recovery} proves the
recovery statements.
Section~\ref{sec:questions} presents questions about larger
minimizing sets, separation from the three-element minimizers,
repeated blow-ups, multipartite and projector minima, other
labeling groups, and recovery. It also gives a certified
four-point obstruction to descent by single exchanges.
Appendix~\ref{app:single-exchange-certificate} justifies
the error bounds used in that certificate.

\section{Bipartite reduction, sharp bounds, and candidate labelings}
\label{sec:bipartite}

\subsection{The bipartite formula and its sharp lower bound}

To calculate the Fourier transform, we first write the edge
indicator as a sum of products of functions of one variable. A pair is an edge exactly when one endpoint
belongs to $C$ and the other does not. Therefore
\begin{equation}\label{eq:indicator-identity}
f_C(x,y)=1_C(x)+1_C(y)-21_C(x)1_C(y).
\end{equation}
Taking its Fourier transform will express the answer in terms of
$S_C$, and summing absolute values will give $D_N(C)$.

\begin{theorem}\label{thm:reduction}
Let $1\leq m\leq n$, let $N=m+n$, and let $C\subset\mathbb Z_N$
have $m$ elements. Then
$$
\operatorname{FR}(f_C)=
\frac{2mn+2(n-m)D_N(C)+2D_N(C)^2}{N\sqrt{2mn}}.
$$
Consequently,
$$
\operatorname{FR}_{\min}(K_{m,n})=
\frac{2mn+2(n-m)d_{m,N}+2d_{m,N}^2}{N\sqrt{2mn}}.
$$
\end{theorem}

\begin{proof}
Taking the Fourier transform of the equation~\eqref{eq:indicator-identity} gives
$$
\widehat f_C(a,b)=S_C(a)1_{\{0\}}(b)+S_C(b)1_{\{0\}}(a)
-\frac{2}{N}S_C(a)S_C(b).
$$
At $(0,0)$ this equals $\frac{2mn}{N}$. If $t\neq0$, the two axis terms are
$$
\widehat f_C(t,0)=\widehat f_C(0,t)=\frac{n-m}{N}S_C(t).
$$
If both frequencies are nonzero, then
$$
|\widehat f_C(a,b)|=\frac{2}{N}|S_C(a)||S_C(b)|.
$$
Summing these three contributions yields
$$
\|\widehat f_C\|_1=
\frac{2mn+2(n-m)D_N(C)+2D_N(C)^2}{N}.
$$
There are $2mn$ ordered edges, so Parseval gives
$\|\widehat f_C\|_2=\sqrt{2mn}$. This proves the first formula.

For $0\leq u<v$, the difference between the numerators at $v$ and
$u$ is
$$
2(v-u)(n-m+v+u)>0.
$$
Thus the formula increases strictly on $[0,\infty)$, including when
$m=n$. Taking the minimum over $C$ proves the second assertion.
\end{proof}

\begin{lemma}\label{lem:affine}
If $u$ is a unit modulo $N$ and $a\in\mathbb Z_N$, then
$$
D_N(uC+a)=D_N(C),
\operatorname{FR}(f_{uC+a})=\operatorname{FR}(f_C).
$$
\end{lemma}

\begin{proof}
The identity
$$
S_{uC+a}(t)=e^{-\frac{2\pi ita}{N}}S_C(ut)
$$
preserves absolute values up to a permutation of the nonzero
frequencies. The second assertion follows from
Theorem~\ref{thm:reduction}.
\end{proof}

The formula reduces our problem to finding the smallest possible
value of $D_N(C)$. The following theorem is the finite cyclic
instance of the classical characterization of norm-one idempotents
in the Fourier algebra. Indeed, for $|C|=m$, the Wiener norm is
\[
\|1_C\|_{A(\mathbb Z_N)}
=\frac1N\sum_{t\in\mathbb Z_N}|S_C(t)|
=\frac{m+D_N(C)}{N}.
\]
The indicator $1_C$ is idempotent under pointwise multiplication,
since $1_C^2=1_C$. For nonempty $C$, its Wiener norm equals one
precisely when $C$ is a coset of a subgroup. This characterization
is included in the stronger small-norm result of Saeki
\cite{Saeki1968}; see its formulation in
\cite[Theorem~1.1(i)]{MudgePham2016}.

We give a self-contained proof in our notation. Fourier inversion
at a point of $C$ gives the lower bound, and the equality case
in the triangle inequality determines when it is attained.

\begin{theorem}\label{thm:coset}
For every nonempty proper subset $C\subset\mathbb Z_N$ of size $m$,
$$
D_N(C)\geq N-m.
$$
Equality holds if and only if $C$ is a coset of a subgroup of
$\mathbb Z_N$.
\end{theorem}

\begin{proof}
Write $\omega=e^{\frac{2\pi i}{N}}$. Fourier inversion gives
$$
1_C(x)=\frac{1}{N}\sum_{t=0}^{N-1}S_C(t)\omega^{tx}.
$$
For $x\in C$, subtract the term $S_C(0)=m$ to obtain
\begin{equation}\label{eq:inversion-at-C}
N-m=\sum_{t=1}^{N-1}S_C(t)\omega^{tx}.
\end{equation}
Taking absolute values proves $D_N(C)\geq N-m$.

Suppose equality holds. For each $x\in C$, equality in the triangle
inequality in the equation~\eqref{eq:inversion-at-C} implies that every nonzero
number $S_C(t)\omega^{tx}$ is positive real. Hence, if $x,y\in C$
and $S_C(t)\neq0$, then
$$
\omega^{t(x-y)}=1.
$$
Define
$$
T=\{t\neq0:S_C(t)\neq0\},
H=\{h\in\mathbb Z_N:\omega^{th}=1\text{ for every }t\in T\}.
$$
The set $H$ is a subgroup. Fix $x_0\in C$. The preceding identity
shows that $C\subset x_0+H$.

The same definition also gives translation invariance. If $h\in H$,
then for every $z\in\mathbb Z_N$,
$$
1_C(z+h)=\frac{1}{N}\sum_{t=0}^{N-1}S_C(t)\omega^{t(z+h)}
=\frac{1}{N}\sum_{t=0}^{N-1}S_C(t)\omega^{tz}=1_C(z).
$$
Only frequencies in $T\cup\{0\}$ contribute. Since $x_0\in C$,
this proves $x_0+H\subset C$. Thus $C=x_0+H$.

Conversely, suppose $C=a+H$ and $|H|=m$. Then $m\mid N$, and the
unique subgroup of this order is
$$
H=\{0,\frac{N}{m},\ldots,(m-1)\frac{N}{m}\}.
$$
A geometric sum gives $|S_C(t)|=m$ when $m\mid t$ and zero
otherwise. Among the nonzero frequencies there are $\frac{N}{m}-1$ multiples
of $m$. Therefore $D_N(C)=m(\frac{N}{m}-1)=N-m$.
\end{proof}

\begin{theorem}\label{thm:sharp}
Let $1\leq m\leq n$. Then
$$
\operatorname{FR}_{\min}(K_{m,n})\geq
\frac{4n^2}{(m+n)\sqrt{2mn}}.
$$
Equality holds if and only if $m\mid n$. When equality holds, a
labeling attains the minimum if and only if the labels on the part of
size $m$ form a coset of the subgroup of order $m$ in $\mathbb Z_{m+n}$.
\end{theorem}

\begin{proof}
Theorem~\ref{thm:coset} gives $D_N(C)\geq n$. Substitute this into
Theorem~\ref{thm:reduction}. At $D_N(C)=n$ the numerator is $4n^2$.
Strict monotonicity shows that equality in the sharp bipartite
bound requires $D_N(C)=n$. By Theorem~\ref{thm:coset}, this is equivalent to the
coset condition. Such a coset has size $m$ exactly when $m\mid N$,
which is equivalent to $m\mid n$.
\end{proof}

\begin{corollary}\label{cor:stars}
Every labeling of $K_{1,n}$ has Fourier ratio
$$
\operatorname{FR}(f_\pi)=\frac{4n^2}{(n+1)\sqrt{2n}}.
$$
Also $\operatorname{FR}_{\min}(K_{n,n})=\sqrt2$, with the minimum
attained precisely when the two parts receive the even and odd labels.
\end{corollary}

\begin{proof}
A singleton is a coset of the trivial subgroup. In $\mathbb Z_{2n}$,
the subgroup of order $n$ consists of the even labels, and its other
coset consists of the odd labels. Apply Theorem~\ref{thm:sharp}.
\end{proof}

\begin{remark}
The adjacency eigenvalues of $K_{m,n}$ are $\sqrt{mn}$,
$-\sqrt{mn}$, and $m+n-2$ zeros. Thus the energy bound gives only
$\sqrt2$. If $r=\frac{n}{m}>1$, our lower bound divided by $\sqrt2$ is
$$
\frac{2r^{\frac{3}{2}}}{1+r}>1,
$$
because $2r^{\frac{3}{2}}>2r>1+r$. The improvement is strict for every
unbalanced complete bipartite graph.
\end{remark}

\subsection{Interval, balanced, and arithmetic progression labelings}
\label{sec:labelings}

For a part of size $m$, we might first try placing the labels next
to one another, or spacing them as evenly as possible around the
cyclic group. These choices give the interval
$$
I_m=\{0,1,\ldots,m-1\}
$$
and the balanced set
$$
B_{m,N}=\{\lfloor\frac{jN}{m}\rfloor:
0\leq j\leq m-1\}.
$$
When $m$ and $N$ are relatively prime, these candidates have exactly
the same Fourier ratio. This fact does not imply that either one
minimizes it.

\begin{proposition}\label{prop:balanced}
If $\gcd(m,N)=1$, then
$$
B_{m,N}=(-m^{-1})I_m\pmod N.
$$
In particular, interval and balanced labelings have the same Fourier
ratio in the coprime case.
\end{proposition}

\begin{proof}
The assertion is immediate for $m=1$. Otherwise, write
$$
b_j=\lfloor\frac{jN}{m}\rfloor,
  r_j=jN-mb_j.
$$
Each $r_j$ lies in $\{0,\ldots,m-1\}$. Since multiplication by $N$
permutes residues modulo $m$, these remainders run through that
entire set. The identity $mb_j=jN-r_j$ gives
$b_j\equiv-m^{-1}r_j\pmod N$. Apply Lemma~\ref{lem:affine}.
\end{proof}

The interval has the explicit value
\begin{equation}\label{eq:interval-D}
D_N(I_m)=\sum_{t=1}^{N-1}
\frac{|\sin(\frac{\pi mt}{N})|}{\sin(\frac{\pi t}{N})}.
\end{equation}
Indeed, this is the absolute value of the usual finite geometric
sum, and the denominator is positive for $1\leq t<N$. A third
candidate, called the parity-tail labeling here, is
$$
P_{m,N}=\{0,q,2q,\ldots,(m-1)q\},
  q=\lfloor\frac{N}{m}\rfloor.
$$
Its elements are distinct because $(m-1)q<N$. The name refers to
the case $q=2$. For general $q$ it is an arithmetic progression with
step $q$.

\begin{proposition}\label{prop:interval-counterexample}
Neither the interval nor the balanced labeling minimizes the Fourier
ratio of $K_{7,9}$.
\end{proposition}

\begin{proof}
Here $N=16$. Let $H=\{0,2,4,\ldots,14\}$ and take
$C=H\setminus\{14\}=P_{7,16}$. The Fourier sum of $1_H$ is $8$ at
$t=0,8$ and zero elsewhere. Subtracting the contribution of $14$
shows that $|S_C(8)|=7$ and $|S_C(t)|=1$ at the other fourteen
nonzero frequencies. Thus
$$
D_{16}(C)=21.
$$
For $I_7$, the seven nonzero even frequencies in the equation
\eqref{eq:interval-D} each contribute $1$. At odd frequencies the
ratio is $|\cot(\frac{\pi t}{16})|$. Pairing terms gives
$$
\begin{aligned}
D_{16}(I_7)
&=7+2\sum_{j\in\{1,3,5,7\}}\cot(\frac{\pi j}{16})\\
&=7+4(\csc(\frac{\pi}{8})
+\csc(\frac{3\pi}{8}))\\
&=7+4\sqrt{8+4\sqrt2}>21.
\end{aligned}
$$
The middle identity uses $\cot x+\tan x=2\csc(2x)$. Squaring the
sum of the two cosecants gives $8+4\sqrt2$,
which is greater than $\frac{49}{4}$. Strict monotonicity in
Theorem~\ref{thm:reduction} proves the assertion for the interval.
Since $\gcd(7,16)=1$, Proposition~\ref{prop:balanced} proves it for
the balanced labeling as well.
\end{proof}

\begin{remark}
We have shown that $P_{7,16}$ is better than the interval and
balanced labelings, without determining whether it is optimal. Even
when $m$ and $N$ are relatively prime, the agreement of the interval
and balanced labelings leaves the minimization problem open.
\end{remark}

\section{All labelings of \texorpdfstring{$K_{2,n}$}{K2,n}}
\label{sec:two}

Assume $n\geq2$ and put $N=n+2$. After translation, the two labels
may be written as $0,d$, where $1\leq d\leq N-1$. Replacing $d$ by
$N-d$ has no effect. We retain all nonzero choices of $d$ and put
$$
g=\gcd(d,N), \qquad k=\frac{N}{g}.
$$
The integer $k\geq2$ is the additive order of $d$ in $\mathbb Z_N$.

For $C=\{0,d\}$, the Wiener norm introduced above satisfies
\[
\|1_C\|_{A(\mathbb Z_N)}=\frac{D_N(C)+2}{N}.
\]
When $k>2$, the cotangent and cosecant evaluations in
\eqref{eq:two-D} are equivalent to the classical formulas for
the norm of the indicator of a union of two cosets, recorded in
\cite[equation~(1.1)]{MudgePham2016}, applied here to the two
singleton cosets $\{0\}$ and $\{d\}$ of the trivial subgroup.
When $k=2$, the set $C$ is itself a subgroup, so its norm is one
by Theorem~\ref{thm:coset}. We give a direct derivation below
and compare the resulting values to determine the minimizing
labelings of $K_{2,n}$.

\begin{theorem}\label{thm:two-values}
For $C=\{0,d\}$,
\begin{equation}\label{eq:two-D}
D_N(C)=
\begin{cases}
2g\cot(\frac{\pi}{2k})-2,&k\text{ even},\\
2g\csc(\frac{\pi}{2k})-2,&k\text{ odd}.
\end{cases}
\end{equation}
For fixed $N$, this quantity decreases strictly with $g$ among
divisors for which $\frac{N}{g}$ is even, and increases strictly with $g$
among divisors for which $\frac{N}{g}$ is odd. Every value from the even
branch is smaller than every value from the odd branch. These
statements also hold for the corresponding Fourier ratios.
\end{theorem}

\begin{proof}
Since $|1+e^{i\theta}|=2|\cos(\frac{\theta}{2})|$,
$$
D_N(C)+2=2\sum_{t=0}^{N-1}
|\cos(\frac{\pi td}{N})|.
$$
Write $d=gd'$ with $\gcd(d',k)=1$. Multiplication by $d'$ permutes
the residues modulo $k$, and absolute cosine has period $k$ in this
sum. Therefore
\begin{equation}\label{eq:cosine-reduction}
D_N(C)+2=2g\sum_{t=0}^{k-1}
|\cos(\frac{\pi t}{k})|.
\end{equation}

If $k=2h$, pairing the terms at $t$ and $k-t$ gives
$$
\sum_{t=0}^{k-1}|\cos(\frac{\pi t}{k})|
=1+2\sum_{t=1}^{h-1}\cos(\frac{\pi t}{k})
=\cot(\frac{\pi}{2k}).
$$
If $k=2h+1$, the analogous calculation gives
$$
\sum_{t=0}^{k-1}|\cos(\frac{\pi t}{k})|
=1+2\sum_{t=1}^{h}\cos(\frac{\pi t}{k})
=\csc(\frac{\pi}{2k}).
$$
Both evaluations follow from the geometric-sum identity
$$
1+2\sum_{t=1}^h\cos(t\theta)
=\frac{\sin((h+\frac{1}{2})\theta)}{\sin(\frac{\theta}{2})}.
$$
Substitution into the equation~\eqref{eq:cosine-reduction} proves the equation
\eqref{eq:two-D}.

To compare the values, set $x=\frac{\pi g}{2N}$. The two expressions for
$D_N(C)+2$ become
$$
\frac{4N}{\pi}x\cot x,
\frac{4N}{\pi}x\csc x.
$$
On $0<x\leq\frac{\pi}{4}$,
$$
\frac{d}{dx}(x\cot x)=\frac{\sin x\cos x-x}{\sin^2x}<0,
\frac{d}{dx}(x\csc x)=\frac{\sin x-x\cos x}{\sin^2x}>0.
$$
The first sign follows from $\sin(2x)<2x$. For the second, the numerator
vanishes at zero and has derivative $x\sin x>0$. Finally,
$x\cot x<1<x\csc x$. This separates the two branches by the value
$\frac{4N}{\pi}-2$. Theorem~\ref{thm:reduction} transfers all comparisons
to the Fourier ratios.
\end{proof}

\begin{corollary}\label{cor:two-count}
Exactly $\tau(N)-1$ different Fourier ratios occur among the labelings
of $K_{2,n}$.
\end{corollary}

\begin{proof}
Every proper positive divisor $g$ of $N$ occurs as $\gcd(d,N)$ by
taking $d=g$. Distinct divisors on the same branch give distinct
values by strict monotonicity. Values on different branches are
separated by Theorem~\ref{thm:two-values}.
\end{proof}

\begin{theorem}\label{thm:two-min}
If $n\geq2$ is even, then
$$
\operatorname{FR}_{\min}(K_{2,n})=
\frac{4n^2}{(n+2)\sqrt{4n}},
$$
and the minimizing pairs are precisely the antipodal pairs
$\{a,a+\frac{N}{2}\}$. If $n\geq3$ is odd, put
$$
\delta_N=2\csc(\frac{\pi}{2N})-2.
$$
Then
$$
\operatorname{FR}_{\min}(K_{2,n})=
\frac{4n+2(n-2)\delta_N+2\delta_N^2}{(n+2)\sqrt{4n}},
$$
and the minimizing pairs are precisely those whose difference is a
unit modulo $N$.
\end{theorem}

\begin{proof}
When $N$ is even, the even branch has its minimum at the largest
possible divisor $g=\frac{N}{2}$, for which $k=2$ and $D_N(C)=N-2=n$.
The odd branch has larger values. When $N$ is odd, only the odd
branch occurs, and its minimum is at $g=1$. Apply
Theorem~\ref{thm:reduction} in both cases.
\end{proof}

\begin{remark}
Parity of $g$ alone does not give the correct monotonicity. For
example, when $N=12$, the even divisors $g=2,4,6$ produce $k=6,3,2$,
respectively. Their $D$ values are
$$
6+4\sqrt3, 14, 10.
$$
They are not decreasing. If $n$ is odd, then $N$ and every divisor
$k$ are odd, so increasing monotonicity in $g$ does hold without
further qualification.
\end{remark}

\section{Three-element minimizing sets}
\label{sec:three}

We now determine the three-element sets $C$ that minimize $D_N(C)$
for every $N\geq4$. If $3\mid N$, Theorem~\ref{thm:coset}
already gives the answer: the minimizers are precisely the
cosets of the subgroup of order $3$. We therefore concentrate
on $3\nmid N$. In this case, we will prove that the minimizers
are precisely
\[
C=t+uC_0=\{t,t+u,t+2u\},
\qquad
C_0=\{0,1,2\},\quad
t\in\mathbb Z_N,\quad u\in\mathbb Z_N^\times.
\]
The complete statement is
Theorem~\ref{thm:three-minimizers-every-order}.
The invertibility of $u$ matters: at composite orders $N$,
an arithmetic progression with three distinct elements need not have
an invertible common difference, and this can increase $D_{N}$, as we shall see.

Proposition~\ref{prop:interval-counterexample} shows that
intervals need not minimize the Fourier ratio for larger
parts. For three labels, however, the interval $C_0$ is
a minimizer whenever $3\nmid N$. To prove this, we must
exclude sets that are not arithmetic progressions and
compare progressions whose common differences have
different additive orders.

The absolute values in $D_N(C)$ make a direct comparison
difficult. Even powers of $|S_C(t)|$ are easier to sum.
After expansion, Fourier orthogonality turns them into
counts of solutions to additive equations in $C$.
An even polynomial lying below $|x|$ therefore gives
a lower bound for $D_N(C)$ in terms of these counts.

This approach follows the MathOverflow argument of
mathworker21, developed in collaboration with Fedja
\cite{Mathworker21Fedja2022}, which proves the three-term
case of the strong Littlewood conjecture on the circle
by combining a degree-six polynomial minorant with
fourth and sixth moment counts. Here we adapt the
counts to additive relations modulo $N$ and construct
rational minorants of degrees six and eight. The
degree-six bound handles sets with sufficiently large
sixth moment; when that moment is smallest, an additional
eighth moment estimate and the degree-eight minorant
give the stronger bound needed.

We first establish the general relation between polynomial
minorants and additive energies. We then calculate the
energies needed for three-element sets and construct the
two minorants. After obtaining the resulting lower bounds,
we compare them with the finite Fourier sum for $C_0$
and compare arithmetic progressions by the additive order
of their common difference. At even orders, sets containing
an antipodal pair require a separate estimate. These
comparisons, together with elementary arguments at the
smallest orders, give the complete classification.

\subsection{Fourier moments and polynomial minorants}

The argument in this subsection applies to every positive
integer $N$ and every nonempty set $C\subset\mathbb Z_N$;
we do not assume that $C$ has three elements or that
$3\nmid N$. Write $m=|C|$. Recall that
\[
S_C(t)=\sum_{c\in C}e^{-2\pi itc/N},
\qquad
D_N(C)=\sum_{t=1}^{N-1}|S_C(t)|.
\]
In particular, $S_C$ is the unnormalized Fourier sum,
and $S_C(0)=m$.

For each positive integer $j$, define
\[
\begin{aligned}
E_j(C)=\bigl|\bigl\{&(x_1,\ldots,x_j,y_1,\ldots,y_j)
\in C^{2j}:\\
&x_1+\cdots+x_j=y_1+\cdots+y_j
\text{ in }\mathbb Z_N\bigr\}\bigr|.
\end{aligned}
\]
Thus $E_j(C)$ counts pairs of ordered $j$-tuples with the
same sum modulo $N$. Repeated entries are allowed.
These quantities are called additive energies. In
particular, $E_1(C)=m$, while
\[
E_2(C)=|\{(a,b,c,d)\in C^4:a+b=c+d
\text{ in }\mathbb Z_N\}|.
\]
For background on additive energy and its relation to
Fourier moments, see \cite[Chapter~4]{TaoVu2006}.

\begin{theorem}\label{thm:polynomial minorant implies estimate}
Let $N$ be a positive integer and let $C\subset\mathbb Z_N$
have $m\geq1$ elements. For every positive integer $j$,
\begin{equation}\label{eq:three-moment-identity}
\sum_{t\in\mathbb Z_N}|S_C(t)|^{2j}=NE_j(C).
\end{equation}
Suppose that, for a positive integer $r$, the polynomial
\[
p(x)=\sum_{j=1}^{r}a_jx^{2j},
\qquad a_j\in\mathbb R,
\]
satisfies $p(x)\leq|x|$ for $|x|\leq m$. Then
\begin{equation}\label{eq:three-general-minorant}
D_N(C)\geq N\sum_{j=1}^{r}a_jE_j(C)-p(m).
\end{equation}
\end{theorem}

\begin{proof}
The orthogonality identity we need is
\[
\sum_{t=0}^{N-1}e^{-2\pi itu/N}
=\begin{cases}
N,&u\equiv0\pmod N,\\
0,&u\not\equiv0\pmod N.
\end{cases}
\]
For the second moment, expanding the square using the $|z|^{2}=z\overline{z}$ identity gives
\[
\begin{aligned}
\sum_{t=0}^{N-1}|S_C(t)|^2
&=\sum_{c,d\in C}\sum_{t=0}^{N-1}
e^{-2\pi it(c-d)/N}\\
&=Nm.
\end{aligned}
\]
Indeed, only the pairs with $c=d$ contribute, and each
contributes $N$.

For the fourth moment, the expansion using $|z|^{4}=z^{2}\overline{z}^{2}$ is
\[
|S_C(t)|^4
=\sum_{a,b,c,d\in C}e^{-2\pi it(a+b-c-d)/N}.
\]
Summing over $t$ leaves precisely the quadruples with
$a+b=c+d$ in $\mathbb Z_N$, each contributing $N$.
Consequently,
\[
\sum_{t=0}^{N-1}|S_C(t)|^4=NE_2(C).
\]

The general identity follows in the same way. We have
\[
\begin{aligned}
|S_C(t)|^{2j}
&=S_C(t)^j\overline{S_C(t)}^{\,j}\\
&=\sum_{x_1,\ldots,x_j\in C}
\sum_{y_1,\ldots,y_j\in C}
e^{-\frac{2\pi it}{N}
(x_1+\cdots+x_j-y_1-\cdots-y_j)}.
\end{aligned}
\]
After summing over $t$, orthogonality assigns weight $N$
to each tuple counted by $E_j(C)$ and weight zero to
every other tuple. This proves
\eqref{eq:three-moment-identity}.

The triangle inequality gives $|S_C(t)|\leq m$, so the
minorant inequality applies at every frequency:
\[
|S_C(t)|\geq p(|S_C(t)|).
\]
Since $D_N(C)$ omits the zero frequency, we sum over
$t\neq0$. Using $S_C(0)=m$ and the moment identity,
we obtain
\[
\begin{aligned}
D_N(C)
&\geq\sum_{t=1}^{N-1}p(|S_C(t)|)\\
&=\sum_{t=0}^{N-1}p(|S_C(t)|)-p(m)\\
&=\sum_{j=1}^{r}a_j
\sum_{t=0}^{N-1}|S_C(t)|^{2j}-p(m)\\
&=N\sum_{j=1}^{r}a_jE_j(C)-p(m).
\end{aligned}
\]
The correction is $p(m)$ because this is the value of
the polynomial at the omitted zero-frequency magnitude.
\end{proof}

For a three-element set, a degree-six minorant
$p(x)=ax^2+bx^4+cx^6$ with $p(3)=3$ therefore gives
\[
D_N(C)\geq\bigl(3a+bE_2(C)+cE_3(C)\bigr)N-3.
\]
A degree-eight minorant uses the same theorem and also
involves $E_4(C)$. There is no sign restriction on the
coefficients in the theorem: the moment identities are
exact. When we subsequently estimate an energy, however,
the sign of its coefficient matters. A positive
coefficient allows us to use a lower bound for that
energy, whereas a negative coefficient requires an
upper bound. We now turn to the additive relations
that provide these estimates for three-element sets.

\subsection{Additive relations in three-element sets}

We now calculate $E_2(C)$ and bound $E_3(C)$ for a
three-element set. The calculation of $E_2(C)$ applies
to both odd and even $N$. When $N$ is even, we call two
elements of $\mathbb Z_N$ antipodal if their difference
is $N/2$.

\begin{lemma}\label{lem:three-energy-all-orders}
Let $N\geq3$ and let $C\subset\mathbb Z_N$ have three elements.
Let $P(C)$ be the number of unordered pairs of distinct
elements $x,y\in C$ satisfying $2x=2y$, and let $M(C)$
be the number of elements $x\in C$ satisfying
\[
2x=y+z,\qquad C=\{x,y,z\}.
\]
Then
\begin{equation}\label{eq:three-energy-all-orders}
E_2(C)=15+2P(C)+4M(C).
\end{equation}

The identity \eqref{eq:three-energy-all-orders} holds
for every $N\geq3$. If $N$ is odd, then $P(C)=0$. If $N$ is even, then
$P(C)$ counts the antipodal pairs in $C$, and $P(C)\leq1$.

\smallskip
\noindent
Under the additional hypothesis that $3\nmid N$,
we have $M(C)\leq1$, and the possibilities are
\[
\begin{array}{c|c|c}
\multicolumn{3}{c}{
\text{\emph{Hypothesis: }}
3\nmid N}\\
\hline
C\text{ is an arithmetic progression}
&
C\text{ contains an antipodal pair}
&
E_2(C)
\\ \hline
\text{no}&\text{no}&15\\
\text{no}&\text{yes}&17\\
\text{yes}&\text{no}&19\\
\text{yes}&\text{yes}&21
\end{array}
\]

In the last case, the common difference of the progression
has additive order $4$.
\end{lemma}

The hypothesis $3\nmid N$ is essential for the table.
For example, when $N=3$ and $C=\mathbb Z_3$, every
element is a midpoint of the other two. Thus, $M(C)=3$
and $P(C)=0$, so the general identity gives
$E_2(C)=27$.

\begin{proof}
We count the ordered solutions of $x+y=u+v$ in $C^4$.
There are three solutions with all entries equal. There
are also six ordered pairs $(x,y)$ with $x\neq y$.
For each of these, $(u,v)$ can be either $(x,y)$ or
$(y,x)$. These give twelve more solutions, for a total
of fifteen.

If both sides contain two distinct elements, the two
pairs share an element because $C$ has only three
elements. Cancelling that element shows that the pairs
agree up to order, so the solution has already been
counted.

Any further solution therefore has a repeated entry on
at least one side. A relation $2x=2u$ with $x\neq u$
contributes two solutions, obtained by interchanging
the two sides. A relation $2x=y+z$ among three distinct
elements contributes four solutions, obtained by ordering
$y,z$ and interchanging the two sides. These contributions
are disjoint and account for all remaining solutions.
This proves \eqref{eq:three-energy-all-orders}.

For odd $N$, multiplication by $2$ is invertible in
$\mathbb Z_N$, so $P(C)=0$. For even $N$, distinct
elements satisfy $2x=2y$ precisely when they differ by
$N/2$. Each element has a unique antipode, so a
three-element set contains at most one antipodal pair.

Suppose now that $3\nmid N$. If two distinct elements
$x,y\in C$ were both midpoints, with $z$ the remaining
element, then
\[
2x=y+z,\qquad 2y=x+z.
\]
Subtracting gives $3(x-y)=0$. Since $3$ is invertible
in $\mathbb Z_N$, this would imply $x=y$, a
contradiction. Thus $M(C)\leq1$.

Moreover, $M(C)=1$ precisely when $C$ is an arithmetic
progression. Indeed, a relation $2x=y+z$ among the three
distinct elements gives
\[
C=\{y,y+d,y+2d\},\qquad d=x-y,
\]
and conversely, the middle term of a three-element
arithmetic progression is a midpoint. The four values
in the table now follow from
\eqref{eq:three-energy-all-orders}.

Finally, write a three-element progression as
$C=\{t,t+d,t+2d\}$, where the three displayed elements
are distinct. Two adjacent terms cannot be antipodal,
since $2d=0$ would identify the endpoints. Thus an
antipodal pair must consist of the endpoints, giving
$2d=N/2$. This is equivalent to $d$ having additive
order $4$.
\end{proof}

In particular, if $3\nmid N$ and $C$ is not an arithmetic
progression and contains no antipodal pair, then
$E_2(C)=15$. We will treat sets containing an antipodal
pair separately, using a direct Fourier estimate in the
final subsection.

We next count three-term sums. An unordered three-term
sum is understood as a multiset: repeated elements are
allowed, but their order is disregarded. Different
multisets can have the same sum in $\mathbb Z_N$.
The following argument gives a lower bound for $E_3(C)$,
characterizes equality, and shows that, when $3\nmid N$,
the next possible value is at least $99$.

\begin{theorem}\label{thm:lower bound for E_3}
Let $N\geq3$ and let $C\subset\mathbb Z_N$ have three
elements. Then $E_3(C)\geq93$, with equality precisely
when the ten multisets of size three drawn from $C$
have pairwise distinct sums in $\mathbb Z_N$.
If $3\nmid N$, then either $E_3(C)=93$ or $E_3(C)\geq99$.
\end{theorem}

\begin{proof}
The ten multisets have the following forms, where
letters appearing together denote distinct elements of $C$:
\[
\begin{array}{c|c|c}
\text{Form of the sum}
&\text{Number of multisets}
&\text{Number of orderings}\\ \hline
3x&3&1\\
2x+y&6&3\\
x+y+z&1&6
\end{array}
\]
For a multiset with $w$ orderings, choosing an ordering
independently for each side gives $w^2$ ordered solutions
counted by $E_3(C)$. Pairing orderings of the same multiset
therefore contributes
\[
3\cdot1^2+6\cdot3^2+1\cdot6^2=93.
\]

All remaining solutions come from distinct multisets
with the same sum. If two such multisets have $m_1$
and $m_2$ orderings, then placing the first on the left
and the second on the right gives $m_1m_2$ solutions.
Interchanging the sides gives another $m_1m_2$ solutions.
Thus each unordered pair of distinct multisets with
the same sum contributes an additional $2m_1m_2$.
These contributions are disjoint, so $E_3(C)\geq93$,
with equality precisely when all ten sums are distinct.

Suppose now that $3\nmid N$. Two distinct pure sums
$3x$ and $3y$ cannot coincide, since multiplication by
$3$ is invertible in $\mathbb Z_N$. Thus, whenever two
distinct multisets have the same sum, at least one
has three or six orderings, and $m_1m_2\geq3$.
Any coincidence therefore contributes at least six
solutions beyond the initial ninety-three. Consequently,
either $E_3(C)=93$ or $E_3(C)\geq99$.
\end{proof}

We now return to the case $3\nmid N$. For the
nonprogressions without antipodal pairs, we can apply
the minorant principle with the exact value $E_2(C)=15$
and the lower bound $E_3(C)\geq93$.
We will choose a degree-six polynomial whose sixth-degree
coefficient is positive, so that any increase in $E_3(C)$
improves the resulting bound for $D_N(C)$. The equality
case $E_3(C)=93$ will require the additional eighth-moment
estimate developed after this construction.

\subsection{Constructing and applying the degree-six minorant}

We seek an even polynomial
\[
p(x)=ax^2+bx^4+cx^6
\]
such that $p(3)=3$, $p(x)\leq |x|$ for $|x|\leq3$, and
$c>0$. For the nonprogressions without antipodal pairs,
we have $E_2(C)=15$ and $E_3(C)\geq93$. Thus the
minorant principle gives
\[
D_N(C)\geq(3a+15b+cE_3(C))N-3.
\]
We therefore want $3a+15b+93c$ to be large. The condition
$c>0$ also ensures that any increase in $E_3(C)$ improves
this lower bound.

For $0<s<3$, let
\[
p_s(x)=a(s)x^2+b(s)x^4+c(s)x^6.
\]
We require $p_s$ to agree with $x$ at $3$ and to be
tangent to $x$ at $s$:
\[
p_s(3)=3,\qquad p_s(s)=s,\qquad p_s'(s)=1.
\]
These conditions give the linear system
\[
\begin{aligned}
9a(s)+81b(s)+729c(s)&=3,\\
s^2a(s)+s^4b(s)+s^6c(s)&=s,\\
2sa(s)+4s^3b(s)+6s^5c(s)&=1.
\end{aligned}
\]
Solving for the coefficients gives
\[
\begin{aligned}
a(s)&=\frac{2s^3+12s^2+54s+81}{6s(s+3)^2},\\
b(s)&=-\frac{4s^3+9s^2+18s+27}{6s^3(s+3)^2},\\
c(s)&=\frac{2s+3}{6s^3(s+3)^2}.
\end{aligned}
\]
In particular, $c(s)>0$ throughout $0<s<3$.

To choose $s$, we consider
\[
L(s)=3a(s)+15b(s)+93c(s).
\]
Numerical optimization places a maximum near
$s=\sqrt{7/2}$, suggesting the nearby rational value
$s=15/8$. The numerical calculation motivates this
choice; the minorant inequality will follow from an
exact factorization valid for every $0<s<3$.

\begin{lemma}\label{lem:specific polynomial}
The polynomial
\[
p(x)=\frac{751}{845}x^2
-\frac{72064}{570375}x^4
+\frac{4096}{570375}x^6
\]
satisfies $p(0)=0$, $p(3)=3$, and
\[
p(x)\leq |x|\qquad (|x|\leq3).
\]
\end{lemma}

\begin{proof}
For $0<s<3$, consider
\[
q_s(x)=x-p_s(x).
\]
The interpolation conditions give
\[
q_s(0)=q_s(3)=q_s(s)=0,\qquad q_s'(s)=0.
\]
Since $0$, $s$, and $3$ are distinct, the factor theorem
shows that $q_s$ is divisible by $x(x-3)(x-s)^2$.
Moreover, $c(s)>0$, so $q_s$ has degree six. Consequently,
there are constants $K\neq0$, $u$, and $v$ such that
\[
q_s(x)=Kx(x-3)(x-s)^2(x^2+ux+v).
\]

We determine the quadratic quotient by comparing three
coefficients. Since
\[
q_s(x)=x-a(s)x^2-b(s)x^4-c(s)x^6,
\]
comparison of the coefficients of $x^6$, $x^5$, and $x$
gives, respectively,
\[
K=-c(s),\qquad
K\bigl(u-(2s+3)\bigr)=0,\qquad
-3Ks^2v=1.
\]
It follows that
\[
u=2s+3,\qquad
v=\frac{1}{3c(s)s^2}
=\frac{2s(s+3)^2}{2s+3}.
\]
The factor theorem already guarantees the existence of
the quotient, and these comparisons determine all three
of its coefficients. Thus no further coefficient checks
are needed, and we obtain
\[
\begin{aligned}
x-p_s(x)
={}&c(s)x(3-x)(x-s)^2\\
&\quad\cdot\left(
x^2+(2s+3)x+\frac{2s(s+3)^2}{2s+3}
\right).
\end{aligned}
\]
For $0\leq x\leq3$, the factors outside the parentheses
are nonnegative. The quadratic factor is positive,
since $s>0$ and $x\geq0$. Therefore $p_s(x)\leq x$ on
$[0,3]$. Since $p_s$ is even, this proves
\[
p_s(x)\leq |x|\qquad (|x|\leq3).
\]

Substituting $s=15/8$ into the coefficient formulas gives
the polynomial in the statement. At this value, the
factorization becomes
\[
x-p(x)=
\frac{x(3-x)(8x-15)^2(64x^2+432x+845)}{570375}.
\]
The values $p(0)=0$ and $p(3)=3$ hold by construction.
\end{proof}

We now apply this polynomial, retaining the contribution
from any excess of $E_3(C)$ over $93$. Write
\[
A=\frac{273581}{190125},\qquad
c_6=\frac{4096}{570375},\qquad
\Gamma=A+6c_6=\frac{281773}{190125}.
\]

\begin{theorem}\label{thm:lower bound for D_N non-AP}
Let $N\geq4$ with $3\nmid N$, and let
$C\subset\mathbb Z_N$ have three elements. If $C$ is
not an arithmetic progression and contains no antipodal
pair, then
\begin{equation}\label{eq:three-sixth-moment-refinement}
D_N(C)\geq
\left(A+c_6(E_3(C)-93)\right)N-3.
\end{equation}
In particular, if $E_3(C)\geq99$, then
\[
D_N(C)\geq\Gamma N-3.
\]
\end{theorem}

\begin{proof}
Lemma~\ref{lem:three-energy-all-orders} gives $E_2(C)=15$.
Apply the polynomial of Lemma~\ref{lem:specific polynomial}
in Theorem~\ref{thm:polynomial minorant implies estimate}.
Since $E_1(C)=3$ and $p(3)=3$, we obtain
\[
\begin{aligned}
D_N(C)
&\geq\left(
\frac{2253}{845}
-\frac{72064}{570375}\cdot15
+\frac{4096}{570375}E_3(C)
\right)N-3\\
&=\left(
\frac{2253}{845}
-\frac{72064}{570375}\cdot15
+\frac{4096}{570375}\cdot93
+c_6(E_3(C)-93)
\right)N-3.
\end{aligned}
\]
The constant part is
\[
\begin{aligned}
\frac{2253}{845}
-\frac{72064}{570375}\cdot15
+\frac{4096}{570375}\cdot93
&=\frac{1520775-1080960+380928}{570375}\\
&=\frac{273581}{190125}=A.
\end{aligned}
\]
This proves \eqref{eq:three-sixth-moment-refinement}.
If $E_3(C)\geq99$, then $c_6>0$ gives
\[
D_N(C)\geq(A+6c_6)N-3=\Gamma N-3,
\]
as claimed.
\end{proof}

Since $E_3(C)\geq93$ by
Theorem~\ref{thm:lower bound for E_3}, the weaker bound
$D_N(C)\geq AN-3$ follows immediately under the same
hypotheses. Retaining the excess in
\eqref{eq:three-sixth-moment-refinement} will also be
useful when $N=7$, where we will obtain a stronger
lower bound for $E_3(C)$.

By Theorem~\ref{thm:lower bound for E_3}, the case not
covered by $E_3(C)\geq99$ is exactly $E_3(C)=93$.
Here the ten unordered three-term sums are distinct.
We will use this restriction to bound $E_4(C)$ from
above and apply a degree-eight minorant to improve
the lower bound for $D_N(C)$.

\subsection{The case of minimum sixth moment}

Suppose that $E_3(C)=93$. By
Theorem~\ref{thm:lower bound for E_3}, the ten unordered
three-term sums from $C$ are distinct. This restricts
coincidences between four-term sums and gives an upper
bound for the fourth additive energy,
\[
\begin{aligned}
E_4(C)=\bigl|\bigl\{&(x_1,\ldots,x_4,y_1,\ldots,y_4)
\in C^8:\\
&x_1+\cdots+x_4=y_1+\cdots+y_4
\text{ in }\mathbb Z_N\bigr\}\bigr|.
\end{aligned}
\]
We will use a degree-eight minorant with a negative
coefficient of $x^8$, so that this upper bound improves
the resulting lower bound for $D_N(C)$.

\begin{lemma}\label{lem:three-higher-energy}
Let $N\geq4$ with $3\nmid N$, and let
$C\subset\mathbb Z_N$ have three elements with $E_3(C)=93$.
Then $N\geq13$, $E_2(C)=15$, and
\[
E_4(C)\leq667.
\]
\end{lemma}

\begin{proof}
Write $C=\{a,b,c\}$. Since the ten unordered three-term
sums are distinct, there must be at least ten residues,
so $N\geq10$. The sum of these ten sums is
\[
10(a+b+c),
\]
since each of $a,b,c$ occurs ten times among them.
If $N=10$, the ten sums would exhaust $\mathbb Z_{10}$.
Their total would therefore be both $0$ and
\[
0+1+\cdots+9=45\equiv5\pmod{10},
\]
a contradiction.

If $N=11$, let $h$ be the unique residue omitted by
the ten sums. Since the sum of all residues modulo
$11$ is zero,
\[
10(a+b+c)=-h\pmod{11}.
\]
Thus $h=a+b+c$ in $\mathbb Z_{11}$, although $a+b+c$
is one of the ten sums. This is again a contradiction.
Finally, $N=12$ is excluded by $3\nmid N$. Hence $N\geq13$.

The distinctness of the three-term sums also determines
$E_2(C)$. If two distinct two-term multisets had the
same sum, adjoining the same element of $C$ to both
would give two distinct three-term multisets with the
same sum. Thus only orderings of the same two-term
multiset contribute to $E_2(C)$. There are three
multisets with both entries equal and three with
distinct entries, giving
\[
E_2(C)=3\cdot1^2+3\cdot2^2=15.
\]

We now count the contributions to $E_4(C)$. The unordered
four-term sums have the following types, where letters
appearing together denote distinct elements of $C$:
\[
\begin{array}{c|c|c}
\text{Type}&\text{Number of multisets}
&\text{Number of orderings}\\ \hline
4x&3&1\\
3x+y&6&4\\
2x+2y&3&6\\
2x+y+z&3&12
\end{array}
\]
Pairing orderings of the same multiset contributes
\[
3\cdot1^2+6\cdot4^2+3\cdot6^2+3\cdot12^2=639.
\]
As in the proof of Theorem~\ref{thm:lower bound for E_3},
two distinct multisets with the same sum and with
$m_1,m_2$ orderings contribute an additional $2m_1m_2$.
It remains to bound these additional contributions.

Two distinct four-term multisets with the same sum
cannot share an element. Otherwise, cancelling one
occurrence of that element from each would leave
two distinct three-term multisets with the same sum.
Their supports, meaning the sets of elements that occur
in them, are therefore disjoint. Since $C$ has only
three elements, at least one of the two multisets
must be a pure sum $4x$.

There are at most three unordered pairs of distinct
four-term multisets with equal sums. To see this, group
the multisets according to their sums. Within each
group, their supports are pairwise disjoint, so the
group has at most three members. A group of three
must consist of the three pure sums. It contributes
three pairs and uses all the pure sums, leaving no
other coincidences. If there is no group of three,
each group containing a coincidence has exactly two
members and uses at least one of the three pure sums.
There are again at most three pairs.

The possible pairs and their additional contributions
to $E_4(C)$ are
\[
\begin{array}{c|c}
\text{Relation}&\text{Contribution}\\ \hline
4x=4y&2\\
4x=3y+z&8\\
4x=2y+2z&12
\end{array}
\]
Here again, letters appearing together denote distinct
elements of $C$.

At most one pair can contribute $12$. For each pure
element $x$, there is only one possible relation of
this type, namely $4x=2y+2z$, where $y,z$ are the other
two elements. Thus two such relations would have
different pure elements. After renaming the elements,
they would be
\[
4x=2y+2z,\qquad 4y=2x+2z.
\]
Subtracting gives $6(x-y)=0$. Since $3$ is invertible
in $\mathbb Z_N$, this implies $2(x-y)=0$, and hence
\[
3x=x+2y.
\]
These are distinct three-term multisets with the same
sum, a contradiction.

There are at most three additional pairs, at most one
contributes $12$, and each other pair contributes at
most $8$. Consequently,
\[
E_4(C)\leq639+12+8+8=667.\qedhere
\]
\end{proof}

We now construct a polynomial that uses this upper
bound. Write
\[
q(x)=a_2x^2+a_4x^4+a_6x^6+a_8x^8
\]
and impose
\[
q(5/3)=5/3,\qquad q'(5/3)=1,\qquad
q(3)=3,\qquad q'(3)=1.
\]
These four linear conditions make $q$ tangent to $x$
at $5/3$ and $3$. Solving for the coefficients gives
the polynomial below; in particular, its coefficient
of $x^8$ is negative, as desired.

\begin{lemma}\label{lem:three-degree-eight-minorant}
The polynomial
\begin{equation}\label{eq:three-degree-eight-minorant}
q(x)=
\frac{19521675x^2-4204351x^4+478053x^6-19521x^8}
{18522000}
\end{equation}
satisfies $q(x)\leq|x|$ for every real $x$, and $q(3)=3$.

If $N\geq4$ with $3\nmid N$ and $C\subset\mathbb Z_N$
has three elements with $E_3(C)=93$, then
\begin{equation}\label{eq:three-degree-eight-bound}
D_N(C)\geq BN-3,
\qquad B=\frac{4489697}{3087000}>\frac{93}{64}.
\end{equation}
\end{lemma}

\begin{proof}
The tangency conditions show that $x-q(x)$ has zeros
of multiplicity at least two at $5/3$ and $3$. It also
vanishes at $0$. Thus, by the factor theorem,
\[
18522000(x-q(x))
=x(3x-5)^2(x-3)^2(Ux^3+Vx^2+Wx+Z)
\]
for some constants $U,V,W,Z$. Comparing the coefficients
of $x^8$, $x^7$, $x^6$, and $x$ gives, respectively,
\[
\begin{aligned}
9U&=19521,\\
9V-84U&=0,\\
9W-84V+286U&=-478053,\\
225Z&=18522000.
\end{aligned}
\]
These equations determine the cubic quotient:
\[
U=2169,\qquad V=20244,\qquad
W=66901,\qquad Z=82320.
\]
We therefore have the exact factorization
\begin{equation}\label{eq:three-degree-eight-factorization}
\begin{aligned}
x-q(x)
={}&\frac{x(3x-5)^2(x-3)^2}{18522000}\\
&\quad\cdot
\left(2169x^3+20244x^2+66901x+82320\right).
\end{aligned}
\end{equation}
Every factor on the right is nonnegative for $x\geq0$.
Hence $q(x)\leq x$ there. Since $q$ is even, it follows
that $q(x)\leq|x|$ for every real $x$. The value
$q(3)=3$ holds by construction and also follows from
the factorization.

For the Fourier estimate, Lemma~\ref{lem:three-higher-energy}
gives $E_2(C)=15$ and $E_4(C)\leq667$.
Apply Theorem~\ref{thm:polynomial minorant implies estimate}
to $q$. Since $E_1(C)=3$, $E_3(C)=93$, and $q(3)=3$,
we obtain
\[
\begin{aligned}
D_N(C)+3
&\geq\frac{N}{18522000}
\left(
19521675\cdot3-4204351\cdot15
+478053\cdot93-19521E_4(C)
\right).
\end{aligned}
\]
The coefficient of $E_4(C)$ is negative, so
$E_4(C)\leq667$ implies
\[
-19521E_4(C)\geq-19521\cdot667.
\]
Consequently,
\[
\begin{aligned}
D_N(C)+3
&\geq\frac{N}{18522000}
\left(
19521675\cdot3-4204351\cdot15
+478053\cdot93-19521\cdot667
\right)\\
&=\frac{26938182}{18522000}N
=\frac{4489697}{3087000}N.
\end{aligned}
\]
Finally,
\[
B-\frac{93}{64}
=\frac{31201}{24696000}>0.
\]
This proves \eqref{eq:three-degree-eight-bound}.
\end{proof}

We now have lower bounds for the nonprogressions without
antipodal pairs in both cases, $E_3(C)\geq99$ and
$E_3(C)=93$. To compare these bounds with the value for
$C_0=\{0,1,2\}$, we next calculate the Fourier sum for
three consecutive labels and compare it with the sums
for other arithmetic progressions.

\subsection{Three consecutive labels and other progressions}

We now calculate the value for $C_0=\{0,1,2\}$ and compare
it with the values for other three-element arithmetic
progressions. For each positive integer $k$, write
\[
W(k)=\frac1k\sum_{t=0}^{k-1}
\left|1+2\cos\left(\frac{2\pi t}{k}\right)\right|,
\qquad
\lambda=\frac13+\frac{2\sqrt3}{\pi}.
\]
Since
\[
S_{C_0}(t)
=1+e^{-2\pi it/N}+e^{-4\pi it/N}
=e^{-2\pi it/N}\left(1+2\cos\frac{2\pi t}{N}\right),
\]
the factor outside the parentheses has absolute value one.
Also, $S_{C_0}(0)=3$, so
\[
D_N(C_0)+3=NW(N).
\]
The following lemma gives an exact formula for $W(k)$
and bounds its difference from $\lambda$. No parity
assumption is needed.

\begin{lemma}\label{lem:three-W-two-sided}
For every integer $k\geq5$ with $3\nmid k$,
\begin{equation}\label{eq:three-W-two-sided}
\frac{\pi^2}{12k^2}
<W(k)-\lambda
<\frac{\pi^2}{4k^2}.
\end{equation}
\end{lemma}

\begin{proof}
Put
\[
r=\left\lfloor\frac{k-1}{3}\right\rfloor,
\qquad
q_t=1+2\cos\left(\frac{2\pi t}{k}\right).
\]
Since $3\nmid k$, none of the $q_t$ vanish. Their positive
terms have indices $0,\pm1,\ldots,\pm r$ modulo $k$.
Taking real parts in the sum of the $k$th roots of unity
gives
\[
\sum_{t=0}^{k-1}q_t=k.
\]
Thus the sum of the absolute values is twice the sum
of the positive terms minus $k$. Since $q_0=3$, we obtain
\[
\begin{aligned}
kW(k)
&=2\left(3+2\sum_{t=1}^r q_t\right)-k\\
&=4r-k+2+
4\frac{\sin((2r+1)\pi/k)}{\sin(\pi/k)}.
\end{aligned}
\]
Here we used the geometric-series identity
\[
1+2\sum_{t=1}^r\cos\left(\frac{2\pi t}{k}\right)
=\frac{\sin((2r+1)\pi/k)}{\sin(\pi/k)}.
\]

Define
\[
\sigma=
\begin{cases}
1,&k\equiv1\pmod3,\\
-1,&k\equiv2\pmod3.
\end{cases}
\]
In the first case, $r=(k-1)/3$; in the second,
$r=(k-2)/3$. Substituting these values gives
\begin{equation}\label{eq:three-W-exact}
W(k)=\frac13+\frac{2\sigma}{3k}
+\frac4k
\frac{\sin(2\pi/3+\sigma\pi/(3k))}
{\sin(\pi/k)}.
\end{equation}

To estimate the error, put
\[
x=\frac{\pi}{k},\qquad
\kappa=\frac{2\sqrt3}{\pi},\qquad
R=\frac{x\cos(x/3)}{\sin x},\qquad
T=\frac{\sin(x/3)}{\sin x}-\frac13.
\]
Since $k\geq5$, we have $0<x<2/3$. Expanding the sine
in \eqref{eq:three-W-exact} gives
\begin{equation}\label{eq:three-W-error}
W(k)-\lambda=\kappa(R-1)-\sigma\frac2kT.
\end{equation}
We bound the two terms separately.

For the lower bound on $R-1$, the Taylor inequalities
\[
\cos(x/3)\geq1-\frac{x^2}{18},
\qquad
\sin x\leq x-\frac{x^3}{6}+\frac{x^5}{120}
\]
give
\[
x\cos(x/3)-\sin x
\geq\frac{x^3}{9}-\frac{x^5}{120}>0.
\]
Since $0<\sin x\leq x$ and $x^2<4/9$,
\[
\begin{aligned}
R-1
&=\frac{x\cos(x/3)-\sin x}{\sin x}\\
&\geq x^2\left(\frac19-\frac{x^2}{120}\right)
\geq\frac{29}{270}x^2.
\end{aligned}
\]
For the upper bound, use $\cos(x/3)\leq1$ and
$\sin x\geq x-x^3/6>0$:
\[
R-1
\leq\frac{x}{x-x^3/6}-1
=\frac{x^2}{6-x^2}
\leq\frac9{50}x^2.
\]

Next, the triple-angle identity
$\sin x=3\sin(x/3)-4\sin^3(x/3)$ gives
\[
T=
\frac{4\sin^2(x/3)}
{3\bigl(3-4\sin^2(x/3)\bigr)}>0.
\]
Using $\sin(x/3)\leq x/3$ and $x^2<4/9$, we obtain
\[
T\leq
\frac{4x^2}{81-12x^2}
\leq\frac{12}{227}x^2
<\frac{x^2}{18}.
\]

Finally, $1<\kappa<7/6$ and $2/k\leq2/5$. For either
value of $\sigma$, equation~\eqref{eq:three-W-error}
therefore implies
\[
\begin{aligned}
W(k)-\lambda
&>\left(\frac{29}{270}-\frac1{45}\right)x^2\\
&=\frac{23}{270}x^2
>\frac{x^2}{12},
\end{aligned}
\]
and
\[
\begin{aligned}
W(k)-\lambda
&<\left(\frac76\cdot\frac9{50}+\frac1{45}\right)x^2\\
&=\frac{209}{900}x^2
<\frac{x^2}{4}.
\end{aligned}
\]
Substituting $x=\pi/k$ proves
\eqref{eq:three-W-two-sided}.
\end{proof}

For the comparisons below, we record the elementary bounds
\[
\sqrt3<\frac{26}{15},\qquad
\pi>\frac{157}{50},\qquad
\pi^2<10.
\]
In particular,
\begin{equation}\label{eq:three-W-rational-bound}
\lambda
<\frac13+\frac{520}{471}
=\frac{677}{471}
<\frac{23}{16}.
\end{equation}

The Fourier sum for a three-element progression depends
on the additive order of its common difference. The
two-sided estimates above allow us to compare this order
with $N$.

\begin{lemma}\label{lem:three-progression-comparison-all-orders}
Let $N\geq4$ with $3\nmid N$, and let
$C=\{h,h+d,h+2d\}\subset\mathbb Z_N$ have three distinct
elements. Then
\[
D_N(C)\geq D_N(C_0),
\]
with equality if and only if $\gcd(d,N)=1$.
\end{lemma}

\begin{proof}
By translation invariance from Lemma~\ref{lem:affine},
we may assume that $h=0$. Put
\[
g=\gcd(d,N),\qquad k=\frac Ng.
\]
The integer $k$ is the additive order of $d$. Since the
three elements are distinct, $k\neq1,2$. Also, $k\mid N$
and $3\nmid N$, so $3\nmid k$. In particular, $k\neq3$,
and hence $k\geq4$.

For every frequency $t$,
\[
\left|1+e^{-2\pi itd/N}+e^{-4\pi itd/N}\right|
=\left|1+2\cos\left(\frac{2\pi td}{N}\right)\right|.
\]
Write $d=gd_0$, so $\gcd(d_0,k)=1$. As $t$ runs through
$\mathbb Z_N$, the residue $td_0$ modulo $k$ takes each
value exactly $g$ times. Therefore
\begin{equation}\label{eq:three-progression-order-reduction}
D_N(\{0,d,2d\})+3=NW(k).
\end{equation}
Since $D_N(C_0)+3=NW(N)$, it remains to compare $W(k)$
and $W(N)$.

If $g=1$, then $k=N$, giving equality. Suppose now that
$g>1$. If $k\geq5$, Lemma~\ref{lem:three-W-two-sided}
applies to both $k$ and $N$ and gives
\[
\begin{aligned}
W(k)-\lambda
&>\frac{\pi^2}{12k^2}
=\frac{g^2\pi^2}{12N^2}\\
&\geq\frac{\pi^2}{3N^2}
>\frac{\pi^2}{4N^2}
>W(N)-\lambda.
\end{aligned}
\]
Thus $W(k)>W(N)$. Here we used only $g\geq2$.

It remains to consider $k=4$. Direct calculation gives
\[
W(4)=\frac{3+1+1+1}{4}=\frac32.
\]
Since $g>1$ and $N=gk$, we have $N\geq8$.
By Lemma~\ref{lem:three-W-two-sided},
\eqref{eq:three-W-rational-bound}, and $\pi^2<10$,
\[
W(N)
<\frac{23}{16}+\frac5{2N^2}
\leq\frac{23}{16}+\frac5{128}
=\frac{189}{128}
<\frac32=W(4).
\]
Thus $W(k)>W(N)$ whenever $g>1$. Equation
\eqref{eq:three-progression-order-reduction} proves the
claimed comparison and its equality characterization.
\end{proof}

We have now determined the minimizing arithmetic
progressions: they are precisely those with invertible
common difference. It remains to exclude nonprogressions.
We will treat sets containing an antipodal pair directly
and compare all other sets with $C_0$ using the minorant
bounds established above.

\subsection{Strict comparison and the final classification}

We first treat the sets containing an antipodal pair. The
following estimate also identifies the equality cases, which
are the progressions of order four encountered above.

\begin{lemma}\label{lem:three-antipodal-bound}
Let $N\geq4$ be even, and suppose that
$C\subset\mathbb Z_N$ has three elements and contains an
antipodal pair. Then
\begin{equation}\label{eq:three-antipodal-bound}
D_N(C)\geq\frac{3N}{2}-3.
\end{equation}
Equality holds precisely when $C$ is an arithmetic progression
whose common difference has additive order $4$.
\end{lemma}

\begin{proof}
Write $N=2M$. By translation invariance, we may assume
\[
C=\{0,M,c\},\qquad c\not\equiv0\pmod M.
\]
Then
\[
S_C(t)=1+(-1)^t+e^{-2\pi itc/N}.
\]
At each of the $M$ odd frequencies, $|S_C(t)|=1$.
At an even frequency $t=2u$,
\[
|S_C(2u)|=\left|2+e^{-2\pi iuc/M}\right|.
\]
Since $S_C(0)=3$, it follows that
\[
\begin{aligned}
D_N(C)+3
&=M+\sum_{u=0}^{M-1}
  \left|2+e^{-2\pi iuc/M}\right|\\
&\geq M+\sum_{u=0}^{M-1}
  \operatorname{Re}\left(2+e^{-2\pi iuc/M}\right)\\
&=3M.
\end{aligned}
\]
The last equality uses
\[
\sum_{u=0}^{M-1}e^{-2\pi iuc/M}=0,
\]
which holds because $c\not\equiv0\pmod M$. This proves
\eqref{eq:three-antipodal-bound}.

Every number $2+e^{-2\pi iuc/M}$ has positive real part.
Equality therefore holds precisely when all these numbers
are real. Taking $u=1$, and using $c\not\equiv0\pmod M$,
this is equivalent to
\[
e^{-2\pi ic/M}=-1.
\]
Thus $M$ is even and $c\equiv M/2\pmod M$. In
$\mathbb Z_{2M}$, the two possibilities are $c=M/2$ and
$c=3M/2$, and both give a three-element progression with
step of order $4$. Conversely, either possibility makes
every character value above equal to $1$ or $-1$, so
equality holds.
\end{proof}

We now compare the lower bounds with the value for
$C_0=\{0,1,2\}$. The following table records the four cases
used in the proof. Throughout the table, $3\nmid N$, and
the remaining nonprogressions are those without antipodal
pairs. We use the constants $\Gamma$ and $B$ defined above,
and write $E_3$ for $E_3(C)$.

\begin{center}
\small
\renewcommand{\arraystretch}{1.3}
\begin{tabular}{@{}p{0.37\linewidth}rl@{}}
\toprule
Case & \shortstack[r]{Lower bound for\\$D_N(C)+3$}
     & Comparison with $W(N)$ \\
\midrule
Antipodal pair, $N\geq7$
  & $\frac32N$ & $W(N)<\frac32$ \\
\addlinespace
Remaining nonprogression,\newline $N=7$
  & $\frac{1487}{975}N$
  & $W(7)<\frac32<\frac{1487}{975}$ \\
\addlinespace
Remaining nonprogression,\newline $N\geq8$, $E_3\geq99$
  & $\Gamma N$ & $W(N)<\frac{189}{128}<\Gamma$ \\
\addlinespace
Remaining nonprogression,\newline $E_3=93$
  & $BN$, with $N\geq13$ & $W(N)<\frac{93}{64}<B$ \\
\bottomrule
\end{tabular}
\end{center}

\begin{proposition}\label{prop:three-strict-nonprogressions}
Let $N\geq7$ with $3\nmid N$, and let
$C\subset\mathbb Z_N$ have three elements. If $C$ is not
an arithmetic progression, then
\[
D_N(C)>D_N(\{0,1,2\}).
\]
\end{proposition}

\begin{proof}
Recall that $D_N(C_0)+3=NW(N)$. By
Lemma~\ref{lem:three-W-two-sided},
\eqref{eq:three-W-rational-bound}, and $\pi^2<10$,
\begin{equation}\label{eq:three-W-small-order-upper}
W(N)<\frac{23}{16}+\frac5{2N^2}.
\end{equation}
In particular, for $N\geq7$,
\[
W(N)<\frac{23}{16}+\frac5{98}
=\frac{1167}{784}<\frac32.
\]

\emph{Case 1: $C$ contains an antipodal pair.}
Lemma~\ref{lem:three-antipodal-bound} gives
\[
D_N(C)+3\geq\frac32N>NW(N)=D_N(C_0)+3.
\]

For the remaining cases, suppose that $C$ has no antipodal
pair. Since $C$ is also not an arithmetic progression,
Theorem~\ref{thm:lower bound for D_N non-AP} applies and gives
\[
D_N(C)+3\geq\bigl[A+c_6(E_3(C)-93)\bigr]N.
\]

\emph{Case 2: $N=7$.}
For each $v\in\mathbb Z_7$, let
\[
r(v)=\left|\{(x,y,z)\in C^3:x+y+z=v\}\right|.
\]
Then
\[
\sum_{v\in\mathbb Z_7}r(v)=27,
\qquad
\sum_{v\in\mathbb Z_7}r(v)^2=E_3(C).
\]
By Cauchy--Schwarz,
\[
27^2\leq7E_3(C).
\]
Since $E_3(C)$ is an integer, this implies
$E_3(C)\geq\lceil729/7\rceil=105$. Therefore
\[
D_N(C)+3
\geq(A+12c_6)N
=\frac{1487}{975}N
>\frac32N
>NW(N),
\]
where
\[
\frac{1487}{975}-\frac32=\frac{49}{1950}>0.
\]

We may now assume $N\geq8$. By
Theorem~\ref{thm:lower bound for E_3}, either $E_3(C)\geq99$
or $E_3(C)=93$.

\emph{Case 3: $N\geq8$ and $E_3(C)\geq99$.}
The degree-six bound gives
\[
D_N(C)+3\geq\Gamma N,
\qquad
\Gamma=A+6c_6=\frac{281773}{190125}.
\]
On the other hand, \eqref{eq:three-W-small-order-upper}
implies
\[
W(N)<\frac{23}{16}+\frac5{2N^2}
\leq\frac{23}{16}+\frac5{128}
=\frac{189}{128}<\Gamma.
\]
The last inequality follows from
\[
\Gamma-\frac{189}{128}
=\frac{133319}{24336000}>0.
\]
Thus $D_N(C)+3>NW(N)$.

\emph{Case 4: $E_3(C)=93$.}
Lemma~\ref{lem:three-higher-energy} gives $N\geq13$,
and Lemma~\ref{lem:three-degree-eight-minorant} gives
\[
D_N(C)+3\geq BN,
\qquad
B=\frac{4489697}{3087000}>\frac{93}{64}.
\]
Using \eqref{eq:three-W-small-order-upper}, we obtain
\[
\begin{aligned}
W(N)
&<\frac{23}{16}+\frac5{2N^2}
\leq\frac{23}{16}+\frac5{338}\\
&<\frac{23}{16}+\frac1{64}
=\frac{93}{64}<B.
\end{aligned}
\]
Here $5/338<1/64$ follows from $320<338$.
Again $D_N(C)+3>NW(N)$.

These cases exhaust the nonprogressions, and each gives
$D_N(C)>D_N(C_0)$.
\end{proof}

The progression comparison and the strict inequality for
nonprogressions now give the complete classification when
$3\nmid N$.

\begin{theorem}\label{thm:three-minimizers-every-order}
Let $N\geq4$ with $3\nmid N$, and put
$C_0=\{0,1,2\}\subset\mathbb Z_N$. Then
\[
\min_{\substack{C\subset\mathbb Z_N\\|C|=3}}D_N(C)
=D_N(C_0).
\]
The minimizing sets are precisely
\[
C=t+uC_0=\{t,t+u,t+2u\},
\qquad
t\in\mathbb Z_N,\quad u\in\mathbb Z_N^\times.
\]
These are exactly the three-element partite sets of
optimal labelings of $K_{3,N-3}$.
\end{theorem}

\begin{proof}
For $N=4$, every three-element set is the complement of
a singleton, so all such sets are translates of $C_0$.
For $N=5$, the complements have two elements. Any two
two-element subsets of $\mathbb Z_5$ are related by an
invertible affine map, so the same holds for their
three-element complements. Thus every triple is an
invertible affine image of $C_0$ in these two cases,
and Lemma~\ref{lem:affine} proves the assertion about $D_N$.

We may now assume $N\geq7$, since $N=6$ is excluded by
$3\nmid N$. Proposition~\ref{prop:three-strict-nonprogressions}
shows that every minimizing triple must be an arithmetic
progression. Lemma~\ref{lem:three-progression-comparison-all-orders}
then shows that its common difference must be invertible,
and that every such progression attains $D_N(C_0)$. This
proves the minimum and its equality classification.

For $N\geq7$, the graph-labeling assertion follows from
Theorem~\ref{thm:reduction} with part sizes $3$ and $N-3$.
For $N=4,5$, apply that theorem to the smaller complementary
part, of size $N-3$. At every nonzero frequency,
\[
S_{\mathbb Z_N\setminus C}(t)=-S_C(t),
\]
because the sum of the corresponding character over
$\mathbb Z_N$ is zero. Hence
\[
D_N(\mathbb Z_N\setminus C)=D_N(C).
\]
The same minimization therefore gives the graph-labeling
assertion in these two cases as well.
\end{proof}

Together with Theorem~\ref{thm:coset}, this determines
the minimizing three-element sets, and hence the optimal
labelings of $K_{3,N-3}$, for every $N\geq4$.
If $3\mid N$, the minimizing sets are precisely the cosets
\[
\{t,t+N/3,t+2N/3\}
\]
of the subgroup of order $3$. If $3\nmid N$, they are
precisely the progressions
\[
\{t,t+u,t+2u\},\qquad u\in\mathbb Z_N^\times,
\]
described in Theorem~\ref{thm:three-minimizers-every-order}. Thus, we know the optimal labelings of $K_{3,n}$ for all $n\geq 1$.

\section{Blow-ups and the energy bound}
\label{sec:blowups}

A labeling of $G$ gives a periodic labeling of every uniform
blow-up of $G$, with the same Fourier ratio. The enlarged graph
may also admit better labelings. We first prove the periodic
calculation, then ask when equality can hold in the energy bound.
We then show that a twofold blow-up can decrease edge complexity
by an arbitrarily large factor, and give a certified finite example.

\begin{definition}
For a positive integer $k$, the uniform independent set blow-up
$G^{(k)}$ replaces each vertex of $G$ by an independent set of size
$k$, and replaces each edge by all edges between the corresponding
two sets. Equivalently,
$$
V(G^{(k)})=V(G)\times\{0,\ldots,k-1\},
$$
and $(v,i)$ is adjacent to $(w,j)$ exactly when $v$ is adjacent to
$w$ in $G$.
\end{definition}

\subsection{Periodic lifting}
The following result is the simple-graph specialization of
\cite[Proposition~6.5]{GuptaIosevichWeighted2026}.
We include the proof to fix the periodic labeling and Fourier
normalization used in the converse below.
\begin{theorem}\label{thm:lift}
Let $G$ have $N$ vertices and $s>0$ edges. For every positive integer
$k$,
$$
\frac{\mathcal E(G)}{\sqrt{2s}}
\leq\operatorname{FR}_{\min}(G^{(k)})
\leq\operatorname{FR}_{\min}(G).
$$
If $G$ attains the energy bound, then both inequalities are equalities.
\end{theorem}

\begin{proof}
Given a labeling $\pi$ of $G$, define a labeling of $G^{(k)}$ by
$$
\pi^{(k)}(x+jN)=(\pi(x),j),
0\leq x<N, 0\leq j<k.
$$
The edge indicator is periodic in each variable, with
$$
f_{\pi^{(k)}}(x+jN,y+\ell N)=f_\pi(x,y).
$$
In its Fourier transform, summation over $j$ gives
$$
\sum_{j=0}^{k-1}e^{-\frac{2\pi iaj}{k}}
=\begin{cases}k,&k\mid a,\\0,&k\nmid a.\end{cases}
$$
The same holds for the second variable. Consequently,
\begin{equation}\label{eq:lift-transform}
\widehat f_{\pi^{(k)}}(a,b)=
\begin{cases}
k\widehat f_\pi(u,v),&a=ku, b=kv, u,v\in\mathbb Z_N,\\
0,&\text{otherwise}.
\end{cases}
\end{equation}
Both Fourier norms are multiplied by $k$, so the ratio is unchanged.
Starting with a minimizing labeling proves the upper bound.

The adjacency matrix of the blow-up is $A_G\otimes J_k$, where
$J_k$ is the all-ones matrix. Its nonzero eigenvalues are the nonzero
eigenvalues of $A_G$ multiplied by $k$, and it has $k^2s$ edges.
Thus
$$
\frac{\mathcal E(G^{(k)})}{\sqrt{2|E(G^{(k)})|}}
=\frac{k\mathcal E(G)}{\sqrt{2k^2s}}
=\frac{\mathcal E(G)}{\sqrt{2s}}.
$$
The energy bound gives the first inequality. If $G$ attains that
bound, the two endpoints coincide.
\end{proof}

\begin{corollary}\label{cor:bipartite-lift}
For positive integers $m\leq n$ and $k$,
$$
\operatorname{FR}_{\min}(K_{km,kn})
\leq\operatorname{FR}_{\min}(K_{m,n}).
$$
If $h=\gcd(m,n)$, then
$$
\operatorname{FR}_{\min}(K_{m,n})
\leq\operatorname{FR}_{\min}
(K_{\frac{m}{h},\frac{n}{h}}).
$$
When $m\mid n$, equality holds in the first inequality for every $k$.
\end{corollary}

\begin{proof}
The graph identities follow directly from the definition of a
blow-up. If $m\mid n$, apply Theorem~\ref{thm:sharp} to both graphs.
Its value is unchanged by multiplying both part sizes by $k$.
\end{proof}

\begin{remark}
For $n=ms$, the graph $K_{m,n}$ is the $m$-fold blow-up of $K_{1,s}$.
When $s>1$, the star does not attain the energy bound. Equality of
the edge complexities in this example comes from the sharper
bipartite estimate, not from the last assertion of
Theorem~\ref{thm:lift}.
\end{remark}

\subsection{The converse for the energy bound}

To prove the converse, we need to recover a labeling of $G$ from
an extremizing labeling of its blow-up. The graph-energy equality
criterion from \cite[Theorem~2.1]{GIIST2026} will show that
vertices with the same neighbors form cosets. Those cosets will
allow us to pass back to $G$.

The underlying matrix result characterizes equality between the
nuclear norm and the entrywise $\ell^1$-norm; see
\cite[Appendix~C, Lemma~3]{ZhuHayashiChen2018}.
We recall the graph-energy criterion below with a proof, and
record its consequence for the nullspace of the adjacency matrix.

\begin{lemma}\label{lem:energy-equality}
Let $A$ be a labeled adjacency matrix and let
$$
U_{tx}=\frac{1}{\sqrt N}e^{-\frac{2\pi itx}{N}},  B=UAU^*.
$$
The labeling attains the energy bound if and only if each row and
each column of $B$ contains at most one nonzero entry. In that case
the nullspace of $A$ is invariant under cyclic translation.
\end{lemma}

\begin{proof}
The two-dimensional transform $UAU$ differs from $B$ by reflection
of its columns, so their entrywise norms agree. Let $W=\operatorname{sgn}(B)$,
defined on eigenvectors by replacing positive eigenvalues by $1$,
negative ones by $-1$, and zeros by $0$. It has operator norm at most
one, and
$$
\mathcal E(G)=\operatorname{tr}(WB)
\leq\sum_{a,b}|W_{ba}||B_{ab}|
\leq\sum_{a,b}|B_{ab}|.
$$
Equality requires $|W_{ba}|=1$ wherever $B_{ab}\neq0$. Each row
and column of a contraction has squared Euclidean norm at most one.
Hence no two such entries can share a row or column. Conversely,
under this support condition the singular values of $B$ are exactly
the absolute values of its nonzero entries, giving equality.

Under the support condition, $\ker B$ is spanned by the coordinate
vectors corresponding to zero columns. Cyclic translation is
diagonal in Fourier coordinates and preserves this subspace.
Conjugating by $U$ proves the last assertion.
\end{proof}

\begin{corollary}\label{cor:regular}
A graph attaining the energy bound is regular.
\end{corollary}

\begin{proof}
For an extremizing labeling, $B_{00}=\frac{2s}{N}>0$. Therefore every other
entry of column zero vanishes. This says that $A$ maps the constant
vector to a constant vector, so all vertex degrees agree.
\end{proof}

\begin{lemma}\label{lem:twins}
Suppose a labeled graph on $\mathbb Z_M$ has an adjacency matrix whose nullspace is invariant under cyclic translation. Then the equivalence classes of vertices having
identical neighbor sets are cosets of a subgroup of $\mathbb Z_M$.
\end{lemma}

\begin{proof}
Let $e_x$ be the coordinate vector at $x$. Since the adjacency matrix
is symmetric, vertices $x,y$ have identical neighbor sets exactly
when $A(e_x-e_y)=0$. Translation invariance of the nullspace shows that
this relation is unchanged when the same element is added to both
labels.

Let $H$ be the class of zero. It contains zero. If $h_1,h_2\in H$,
then $h_1$ and $h_2$ are equivalent. Translating by $-h_2$ shows that
$h_1-h_2\in H$. Thus $H$ is a subgroup. Translation invariance gives
$x\sim y$ exactly when $x-y\in H$, proving the assertion.
\end{proof}

\begin{theorem}\label{thm:blowup-converse}
Let $G$ have at least one edge and let $k\geq1$. The graph $G^{(k)}$
attains its energy bound if and only if $G$ attains its energy bound.
\end{theorem}

\begin{proof}
One direction is Theorem~\ref{thm:lift}. For the other, let $G$ have
$N$ vertices and choose an extremizing labeling of $G^{(k)}$ by
$\mathbb Z_{kN}$. By Lemmas~\ref{lem:energy-equality} and
\ref{lem:twins}, the classes of vertices with identical neighbor
sets are cosets of one subgroup $H$.

In the abstract blow-up, two vertices $(v,i)$ and $(w,j)$ have the
same neighbor set exactly when $v$ and $w$ do in $G$. Hence every
such class in $G^{(k)}$ has size $k$ times the size of the
corresponding class in $G$. In particular, $k\mid |H|$.

The cyclic group $\mathbb Z_{kN}$ has a unique subgroup of order $k$,
namely
$$
K=\{0,N,2N,\ldots,(k-1)N\}.
$$
Since $k\mid|H|$, this subgroup lies in $H$. Thus vertices whose
labels differ by $N$ have identical neighbor sets. The adjacency
matrix is periodic in each variable with period $N$ and descends to
an adjacency matrix $A_0$ on $\mathbb Z_N$.

We verify that $A_0$ represents $G$. Vertices with the same neighbor
set form an independent set. If two distinct such vertices were
adjacent, their identical neighbor sets would force a loop. Between
two classes either all edges or no edges are present. Passing from
$G^{(k)}$ to $A_0$ divides the size of each class by $k$ and preserves
these complete or empty relations. Those are exactly the class
sizes and relations in $G$. Therefore $A_0$ is a labeling of a
graph isomorphic to $G$.

The periodic transform calculation in the equation~\eqref{eq:lift-transform} now
gives $\operatorname{FR}(A_0)=\operatorname{FR}(A_{G^{(k)}})$.
The energy quotient is also unchanged by a blow-up, as computed in
Theorem~\ref{thm:lift}. Since the latter labeling attains its energy
bound, $A_0$ attains the energy bound of $G$.
\end{proof}

\begin{corollary}
If $G$ does not attain its energy bound, then for every $k\geq1$,
$$
\operatorname{FR}_{\min}(G^{(k)})>
\frac{\mathcal E(G)}{\sqrt{2|E(G)|}}.
$$
In particular this holds for every nonregular graph with at least
one edge.
\end{corollary}

\begin{proof}
The energy quotient is unchanged by a blow-up. Equality would
contradict Theorem~\ref{thm:blowup-converse}. Apply
Corollary~\ref{cor:regular} for the final statement.
\end{proof}

\subsection{A strict decrease for complete bipartite graphs}

The converse just proved concerns attainment of the energy bound.
The next result shows that the edge complexity itself can decrease
by an arbitrarily large factor under a twofold blow-up.

\begin{theorem}\label{thm:strict-family}
For an odd prime $p$, put
$$
G_p=K_{\frac{p-1}{2},\frac{p+1}{2}}.
$$
There is an absolute constant $c>0$ such that, for sufficiently large
$p$,
$$
\operatorname{FR}_{\min}(G_p)
\geq c\frac{\log p}{(\log\log p)^3}.
$$
On the other hand, for every odd prime $p$,
\begin{equation}\label{eq:strict-upper}
\operatorname{FR}_{\min}(G_p^{(2)})
\leq\sqrt2(5-\frac{1}{p})
\sqrt{\frac{p-1}{p+1}}<5\sqrt2.
\end{equation}
Consequently, for all sufficiently large $p$ the blow-up inequality
is strict, and the ratio
$\frac{\operatorname{FR}_{\min}(G_p)}{\operatorname{FR}_{\min}(G_p^{(2)})}$
tends to infinity.
\end{theorem}

\begin{proof}
To bound the Fourier ratio before the blow-up, we use Sanders's
Theorem~1.4 in \cite{Sanders2007}, which improves the theorem of
Green and Konyagin \cite{GreenKonyagin2009}. In our notation it gives
\begin{equation}\label{eq:sanders-input}
\frac{1}{p}\sum_{t=0}^{p-1}|S_C(t)|
\geq c_0(\frac{\log p}{(\log\log p)^3})^{\frac{1}{2}}
\end{equation}
uniformly over subsets of $\mathbb Z_p$ whose density lies between
$\frac{1}{3}$ and $\frac{2}{3}$, for sufficiently large primes. The transform in that
theorem uses probability measure. The factor $\frac{1}{p}$ in the equation
\eqref{eq:sanders-input} makes the normalization explicit.

Take $|C|=m=\frac{p-1}{2}$ and set $\alpha=\frac{m}{p}$ and $\eta=\frac{D_p(C)}{p}$.
The left side of the equation~\eqref{eq:sanders-input} is $\alpha+\eta$.
Since $\alpha\leq\frac{1}{2}$, we have, after increasing the lower threshold
for $p$,
$$
\eta\geq\frac{c_0}{2}
(\frac{\log p}{(\log\log p)^3})^{\frac{1}{2}}.
$$
Theorem~\ref{thm:reduction}, divided through by $p^2$, reads
\begin{equation}\label{eq:density-formula}
\operatorname{FR}(f_C)=
\frac{2\alpha(1-\alpha)+2(1-2\alpha)\eta+2\eta^2}{\sqrt{2\alpha(1-\alpha)}}.
\end{equation}
Here $\alpha\leq\frac{1}{2}$, so every numerator term is nonnegative.
Also $\sqrt{2\alpha(1-\alpha)}\leq\frac{1}{\sqrt2}$. Keeping the last
term proves the asserted lower bound, uniformly over all labelings.

For the blow-up, which is $K_{p-1,p+1}$ on $\mathbb Z_{2p}$, take
the subgroup of even labels with zero removed,
$$
C=\{2,4,\ldots,2p-2\}.
$$
The full subgroup has Fourier sum $p$ at $t=0,p$ and zero elsewhere.
Thus $S_C(p)=p-1$, while $S_C(t)=-1$ at the other $2p-2$ nonzero
frequencies. It follows that
$$
D_{2p}(C)=3(p-1).
$$
Substituting $m=p-1$, $n=p+1$, and this value into
Theorem~\ref{thm:reduction} gives
$$
\operatorname{FR}(f_C)=
\frac{4(p-1)(5p-1)}{2p\sqrt{2(p^2-1)}}
=\sqrt2(5-\frac{1}{p})\sqrt{\frac{p-1}{p+1}}.
$$
Each of the last two factors is strictly smaller than $5$ and $1$,
respectively. This proves the bound in the equation~\eqref{eq:strict-upper} and the conclusion.
\end{proof}

The improvement comes from the subgroups available for labeling.
Before the blow-up the
ambient cyclic group has prime order, and indicators with density
near one half must have a growing Fourier $\ell^1$-norm. After the
blow-up, a subgroup of index two is available. Removing one point
from that subgroup gives exactly the required part size while
keeping the normalized Fourier norm bounded. The periodic lift is
therefore not always the best labeling of the enlarged graph.

\subsection{A certified finite counterexample}
\label{sec:certificate}

Theorem~\ref{thm:strict-family} gives strict decrease for all
sufficiently large primes, but does not give a small numerical
example. For $p=23$, we can check every labeling. We use rational
approximations with proved error bounds, so the comparison follows
from exact integer inequalities.

\begin{proposition}[Computer-assisted]\label{prop:finite-counterexample}
The twofold blow-up of $K_{11,12}$ has strictly smaller edge
complexity, with
$$
\operatorname{FR}_{\min}(K_{22,24})
<\operatorname{FR}_{\min}(K_{11,12}).
$$
\end{proposition}

\begin{proof}
In $\mathbb Z_{46}$ choose $C=\{2,4,\ldots,44\}$. The calculation
in Theorem~\ref{thm:strict-family} gives $D_{46}(C)=66$. We certify
that every eleven-element subset $A\subset\mathbb Z_{23}$ satisfies
\begin{equation}\label{eq:certified-D}
D_{23}(A)\geq\frac{33524494}{1000000}>33.
\end{equation}
Both graphs have the same smaller-part density $\frac{11}{23}$, while
$$
\frac{D_{46}(C)}{46}=\frac{33}{23}
<\frac{D_{23}(A)}{23}.
$$
Strict monotonicity in the equation~\eqref{eq:density-formula} then proves the
proposition.

To verify the lower bound, translation invariance
allows us to assume $0\in A$, leaving
$\binom{22}{10}=646646$ sets. Put $Q=10^6$. For each residue $j$,
the accompanying program constructs integers $a_j,b_j$ satisfying
$$
|\frac{a_j}{Q}-\cos(\frac{2\pi j}{23})|
<\frac{2}{Q},
|\frac{b_j}{Q}-\sin(\frac{2\pi j}{23})|
<\frac{2}{Q}.
$$
They are obtained from Machin's identity
$$
\pi=16\arctan(\frac{1}{5})
-4\arctan(\frac{1}{239})
$$
using thirty terms of each arctangent series, followed by forty
terms of each sine and cosine series, all evaluated as rational
numbers. The error in the approximation to $\pi$ is bounded by
$$
\frac{16}{61\cdot5^{61}}+\frac{4}{61\cdot239^{61}}.
$$
For arguments of absolute value at most $7$, the degree-$79$
Taylor remainder is at most $\frac{7^{80}}{80!}$. Twice the displayed
error for $\pi$, plus this remainder, is less than $\frac{1}{Q}$.
Taking the integer floor after multiplying by $Q$ accounts for the
other $\frac{1}{Q}$ in the coordinate estimates.

For each set $A$ and $1\leq t\leq11$, form the integer sums
$$
X_t=\sum_{c\in A}a_{tc\bmod23},
Y_t=\sum_{c\in A}b_{tc\bmod23},
R_t=\lfloor\sqrt{X_t^2+Y_t^2}\rfloor.
$$
The complex approximation error is less than
$\frac{22\sqrt2}{Q}<\frac{33}{Q}$. Therefore
$$
|S_A(t)|\geq\frac{R_t-33}{Q}.
$$
Conjugation pairs the nonzero frequencies, so
$$
D_{23}(A)\geq\frac{2}{Q}\sum_{t=1}^{11}(R_t-33).
$$
The minimum integer numerator over the $646646$ sets is $33524494$,
which proves the bound in the equation~\eqref{eq:certified-D}. Each square-root floor is checked
by the exact inequalities
$R_t^2\leq X_t^2+Y_t^2<(R_t+1)^2$. All proof-relevant calculations use arbitrary-precision integers
and exact rational arithmetic. The accompanying computational
supplement contains \nolinkurl{verify_blowup_counterexample.py}
and its shared arithmetic module \nolinkurl{certified_fourier.py},
together with instructions, tests, and recorded results.
\end{proof}

The table gives rounded values from the formulas above. The last
column states whether the value is a proved minimum or comes from
a particular labeling. The strict comparison for the last two
graphs follows from the rational certificate, not from the rounded
values.

\begin{center}
\small
\begin{tabular}{@{}lll@{}}
\toprule
Graph & Fourier ratio & Status\\
\midrule
$K_{1,4}$ & $4.525483$ & Every labeling\\
$K_{2,3}$ & $3.518619$ & Exact minimum, Theorem~\ref{thm:two-min}\\
$K_{2,5}$ & $5.097900$ & Exact minimum, Theorem~\ref{thm:two-min}\\
$K_{7,9}$ & $6.080193$ & Upper bound from $P_{7,16}$\\
$K_{11,12}$ & $6.901001$ & Upper bound from $I_{11}$\\
$K_{22,24}$ & $6.711162$ & Upper bound from the even subgroup minus zero\\
\bottomrule
\end{tabular}
\end{center}

\section{Complete multipartite graphs and joins}
\label{sec:multipartite}

For complete bipartite graphs with two partite sets, the nonzero Fourier sums of the two indicators are negatives of one another. With more than two partite sets, we must keep track of
how their Fourier sums combine more carefully. The calculation still gives an
exact formula, but the separate absolute values no longer determine
the answer.

Let $r\geq2$, let $n_1,\ldots,n_r$ be positive integers, and set
\[
N=\sum_{j=1}^r n_j,\qquad
a=\min_j n_j,\qquad b=\max_j n_j,
\]
\[
\mu=\frac{1}{N}\sum_{j=1}^r n_j^2,\qquad
2s=N^2-\sum_{j=1}^r n_j^2=N(N-\mu).
\]
A labeling of $K_{n_1,\ldots,n_r}$ is determined by a partition
$\mathbb Z_N=C_1\cup\cdots\cup C_r$ with $|C_j|=n_j$. Write
\[
S_j(t)=\sum_{x\in C_j}e^{-2\pi itx/N},
\qquad D_j=\sum_{t\neq0}|S_j(t)|.
\]
We first record the labeling formula. We then prove a sharp lower
bound and classify all part sizes and labelings attaining it.

\subsection{The labeling formula}

\begin{theorem}\label{thm:multipartite-formula}
For the labeling determined by the partition above,
\begin{equation}\label{eq:multipartite-formula}
\operatorname{FR}(f)=\frac{1}{N\sqrt{2s}}
\left(
2s+2\sum_{t\neq0}\left|\sum_{j=1}^r n_jS_j(t)\right|
+\sum_{\alpha,\beta\neq0}
\left|\sum_{j=1}^r S_j(\alpha)S_j(\beta)\right|
\right).
\end{equation}
For every real number $u$,
\begin{equation}\label{eq:multipartite-upper}
\operatorname{FR}(f)\leq
\frac{2s+2\sum_{j=1}^r|n_j-u|D_j+\sum_{j=1}^rD_j^2}
{N\sqrt{2s}}.
\end{equation}
\end{theorem}

\begin{proof}
The edge indicator is
\[
f(x,y)=1-\sum_{j=1}^r1_{C_j}(x)1_{C_j}(y),
\]
so
\[
\widehat f(\alpha,\beta)
=N1_{\{0\}}(\alpha)1_{\{0\}}(\beta)
-\frac1N\sum_{j=1}^rS_j(\alpha)S_j(\beta).
\]
Since $S_j(0)=n_j$, the value at the origin is
\[
\widehat f(0,0)
=N-\frac1N\sum_{j=1}^r n_j^2
=\frac{2s}{N}.
\]
On either nonzero axis the value is
$-N^{-1}\sum_{j=1}^r n_jS_j(t)$, where $t\neq0$.
Off the axes, where $\alpha,\beta\neq0$, the value is
\[
-\frac1N\sum_{j=1}^rS_j(\alpha)S_j(\beta).
\]
Summing absolute values over these three regions and using
$\|\widehat f\|_2=\|f\|_2=\sqrt{2s}$ proves
\eqref{eq:multipartite-formula}.

Since the parts partition $\mathbb Z_N$, we have
$\sum_{j=1}^rS_j(t)=0$ for $t\neq0$. Thus, for every real number $u$,
\[
\sum_{j=1}^r n_jS_j(t)
=\sum_{j=1}^r(n_j-u)S_j(t).
\]
The triangle inequality therefore gives
\[
\sum_{t\neq0}\left|\sum_{j=1}^r n_jS_j(t)\right|
\leq
\sum_{j=1}^r |n_j-u|\sum_{t\neq0}|S_j(t)|
=\sum_{j=1}^r |n_j-u|D_j.
\]
For the off-axis terms, another application of the triangle
inequality gives
\[
\begin{aligned}
\sum_{\alpha,\beta\neq0}
\left|\sum_{j=1}^r S_j(\alpha)S_j(\beta)\right|
&\leq
\sum_{j=1}^r
\left(\sum_{\alpha\neq0}|S_j(\alpha)|\right)
\left(\sum_{\beta\neq0}|S_j(\beta)|\right)\\
&=\sum_{j=1}^r D_j^2.
\end{aligned}
\]
Substituting these two estimates into
\eqref{eq:multipartite-formula} proves
\eqref{eq:multipartite-upper}.
\end{proof}

\begin{remark}
For $r=2$, $S_2(t)=-S_1(t)$ at nonzero frequencies, so
\eqref{eq:multipartite-formula} reduces exactly to the bipartite
formula in Theorem~\ref{thm:reduction}. For more parts, the phases
in the sums over $j$ matter. The numbers $D_j$ alone give the
upper bound \eqref{eq:multipartite-upper}, but need not determine
the exact Fourier ratio.
\end{remark}

\subsection{A sharp lower bound}

Let $J$ denote the all-ones matrix, let $\mathbf1$ be the constant
vector, and let
\[
U_{tx}=\frac1{\sqrt N}e^{-2\pi itx/N}.
\]
The matrices $UAU$ and $UAU^*$ differ by reflection of their columns,
so their entrywise $\ell^1$-norms agree. We write
$\|M\|_{\ell^1}=\sum_{x,y}|M_{xy}|$ for this matrix norm.

\begin{theorem}\label{thm:multipartite-sharp-bound}
For every complete multipartite graph with at least two parts,
\begin{equation}\label{eq:multipartite-sharp-bound}
\operatorname{FR}_{\min}(K_{n_1,\ldots,n_r})
\geq \frac{2(N+\mu-2a)}{\sqrt{2s}}
=\frac{2\left(N^2+\sum_jn_j^2-2aN\right)}
{N\sqrt{N^2-\sum_jn_j^2}}.
\end{equation}
\end{theorem}

\begin{proof}
Fix a labeling and put
\[
S=\sum_{j=1}^r1_{C_j}1_{C_j}^*,\qquad
P=I-\frac1N J,\qquad R=PSP.
\]
Then $A=J-S$. If $n(x)$ is the size of the part containing $x$,
direct expansion gives
\begin{equation}\label{eq:multipartite-centered-kernel}
R(x,y)=1_{\{x,y\text{ belong to the same part}\}}
-\frac{n(x)+n(y)}N+\frac{\mu}{N}.
\end{equation}
The matrix $R$ is positive semidefinite and annihilates $\mathbf1$.

Write $B=UAU^*$, $\beta_t=B_{t0}$ for $t\neq0$, and
\[
T=(URU^*)_{t,u\neq0}.
\]
The corresponding off-axis block of $B$ is $-T$, and
$B_{00}=N-\mu$. Hence
\begin{equation}\label{eq:multipartite-block-split}
\|\widehat A\|_1
=N-\mu+2\sum_{t\neq0}|\beta_t|+\|T\|_{\ell^1}.
\end{equation}
Fourier inversion of the degree function yields
\begin{equation}\label{eq:multipartite-axis-inversion}
\sum_{t\neq0}\beta_t e^{2\pi itx/N}=\mu-n(x).
\end{equation}
At a vertex in a part of size $a$, this gives
\[
\sum_{t\neq0}|\beta_t|\geq\mu-a.
\]
At the same vertex, matrix Fourier inversion and
\eqref{eq:multipartite-centered-kernel} give
\[
\|T\|_{\ell^1}\geq NR(x,x)=N-2a+\mu.
\]
Substitution into \eqref{eq:multipartite-block-split} proves
$\|\widehat A\|_1\geq2(N+\mu-2a)$. Divide by
$\|A\|_2=\sqrt{2s}$ and minimize over labelings.
\end{proof}

The equality condition is more restrictive than the bound itself
suggests. Besides the balanced graphs, it permits exactly one large
part, with all other parts equal and an additional divisibility
condition.

\begin{theorem}\label{thm:multipartite-sharp-equality}
Equality in \eqref{eq:multipartite-sharp-bound} holds if and only if,
after reordering the parts, one of the following occurs:
\begin{enumerate}
\item $n_1=\cdots=n_r=a$;
\item $n_1=\cdots=n_h=a$ and $n_{h+1}=kha$, where
$h=r-1$ and $k\geq1$ is an integer.
\end{enumerate}
The two cases overlap only when $r=2$ and $k=1$.

In the balanced case, the extremizing labelings are exactly those
whose parts are the cosets of the subgroup of order $a$ in
$\mathbb Z_N$. In the unbalanced case, let $H\subseteq L$ be the
subgroups of orders $a$ and $ha$, respectively. A labeling is
extremizing if and only if, for some $z\in\mathbb Z_N$, the small
parts are the $h$ cosets of $H$ contained in $z+L$, and the large
part is $\mathbb Z_N\setminus(z+L)$.
\end{theorem}

\begin{proof}
Suppose first that a labeling attains the lower bound
\eqref{eq:multipartite-sharp-bound}. Translate a vertex
in a smallest part to zero. Use $R$ and $T$ from the preceding
proof. Since
\[
\sum_{t,u\neq0}T_{tu}=NR(0,0)=N-2a+\mu,
\]
equality in the estimate for $\|T\|_{\ell^1}$ forces every entry of
$T$ to be real and nonnegative.

We first determine the possible part sizes. Let $\lambda>0$ be the
largest eigenvalue of $R$. The nonzero eigenvalues of $R$ are those
of $T$. Since $T$ is real symmetric and entrywise nonnegative,
replacing a real unit top eigenvector by its coordinatewise
absolute value cannot decrease its Rayleigh quotient. Thus $T$
has a nonzero eigenvector $z\geq0$ for $\lambda$. Extend $z$ by zero at frequency
zero and put $v=U^*z$. Then
\begin{equation}\label{eq:multipartite-perron-maximum}
v(0)=\frac1{\sqrt N}\sum_{t\neq0}z_t>0,
\qquad |v(x)|\leq v(0)\quad(x\in\mathbb Z_N).
\end{equation}
The matrix $S$ maps every vector to a vector constant on
each part. The projection $P$ subtracts the mean of a vector,
preserving constancy on each part and giving mean zero.
Since $Rv=\lambda v$ with $\lambda>0$, we have
\[
v=\frac1\lambda PSPv.
\]
Thus $v$ is constant on each part and has mean zero.
Writing $v_j$ for its value on $C_j$, the eigenvalue equation is
\begin{equation}\label{eq:multipartite-compressed-eigenvector}
(n_j-\lambda)v_j=\gamma,\qquad
\gamma=\frac1N\sum_i n_i^2v_i.
\end{equation}

Assume the graph is unbalanced. If the largest part size $b$ occurs
at least twice, then $\lambda=b$. Indeed, $S\leq bI$, while a
mean-zero difference of indicators of two largest parts is an
eigenvector of $R$ with eigenvalue $b$. Equality in the Rayleigh
bound shows that every $b$-eigenvector is supported on the largest
parts: it must lie in $\mathbf1^\perp$ and in the $b$-eigenspace
of $S$. It therefore vanishes at zero, contradicting
\eqref{eq:multipartite-perron-maximum}.

The largest part is thus unique. Let $d$ be the second-largest
size. Then
\begin{equation}\label{eq:multipartite-top-eigenvalue}
d<\lambda<b.
\end{equation}
For the lower inequality, choose parts of sizes $d$ and $b$ and use
the mean-zero vector $1_{C_d}/d-1_{C_b}/b$. Its Rayleigh quotient
for $R$ is $2bd/(b+d)>d$. For the upper inequality, equality
$\lambda=b$ would require a mean-zero vector in the one-dimensional
$b$-eigenspace of $S$, which is impossible.

In particular $\gamma\neq0$ in
\eqref{eq:multipartite-compressed-eigenvector}. If some part had
size $a<n_j<b$, then
\[
|v_j|=\frac{|\gamma|}{\lambda-n_j}
>\frac{|\gamma|}{\lambda-a}=v(0),
\]
contradicting \eqref{eq:multipartite-perron-maximum}. Thus there are
$h=r-1$ parts of size $a$ and one part of size $b>a$.

We next obtain the divisibility condition and the labeling
structure, treating the balanced case at the same time.
The matrix $URU^*$ has nonnegative entries: its off-axis block
is $T$, and both zero-frequency axes vanish because
$R\mathbf1=0$ and $R=R^*$. The two-dimensional Fourier
transform $\widehat R=URU$ is obtained from $URU^*$ by
reflecting its columns. Hence
\[
\widehat R(t,u)\geq0
\qquad(t,u\in\mathbb Z_N).
\]

Fourier inversion, with our unitary normalization, gives
\[
R(x,y)=\frac1N\sum_{t,u\in\mathbb Z_N}
\widehat R(t,u)e^{2\pi i(tx+uy)/N}.
\]
Since $R$ is real-valued, taking real parts yields
\[
R(0,0)-R(x,y)
=\frac1N\sum_{t,u\in\mathbb Z_N}
\widehat R(t,u)
\left(1-\cos\frac{2\pi(tx+uy)}N\right).
\]
Every summand is nonnegative. Thus $R(x,y)=R(0,0)$
if and only if
\[
e^{2\pi i(tx+uy)/N}=1
\quad\text{whenever }\widehat R(t,u)>0.
\]
For each $(t,u)$, the character
\[
\chi_{t,u}(x,y)=e^{2\pi i(tx+uy)/N}
\]
is a group homomorphism from $\mathbb Z_N^2$ to the unit
circle under multiplication. Its kernel is therefore a subgroup.
Consequently,
\[
\mathcal M
=\{(x,y)\in\mathbb Z_N^2:R(x,y)=R(0,0)\}
=\bigcap_{\widehat R(t,u)>0}\ker\chi_{t,u}
\]
is a subgroup of $\mathbb Z_N^2$.

Formula \eqref{eq:multipartite-centered-kernel} identifies this set:
\[
\mathcal M=\{(x,y):x,y\text{ belong to the same part of size }a\}.
\]
If $C_0$ is the part containing zero, then
\[
H=\{x:(x,0)\in\mathcal M\}=C_0
\]
is a subgroup of order $a$. Also
\[
L=\{x:(x,x)\in\mathcal M\}
\]
is a subgroup consisting of the union of all smallest parts.
If $x$ lies in a smallest part $C_j$, then $(x,x)\in\mathcal M$.
For every $h_0\in H$, the subgroup property gives
$(x+h_0,x)=(x,x)+(h_0,0)\in\mathcal M$. Thus $x+H\subset C_j$.
Since both sets have size $a$, each smallest part is an $H$-coset.

In the balanced case, these are all cosets of $H$ in $\mathbb Z_N$.
In the unbalanced case, $|L|=ha$, so $ha\mid N=ha+b$ and therefore
$ha\mid b$. The smallest parts are exactly the $H$-cosets in $L$,
and the remaining part is its complement. Undoing the translation
gives the stated necessary labeling conditions.

Conversely, in the balanced case label the parts by the cosets of
the subgroup of order $a$. The matrix $S$, and hence $R$, is
circulant. Its Fourier conjugate is diagonal and nonnegative.
The axes of $UAU^*$ vanish because all degrees are equal.
Thus equality holds in the proof of
Theorem~\ref{thm:multipartite-sharp-bound}.

For the unbalanced case write $b=qa$, where $h\mid q$, and set
$M=h+q$. First work in $\mathbb Z_M$. Let $H_0$ be its subgroup of
order $h$, assign its elements to the $h$ singleton parts, and
assign its complement to the part of size $q$. Its character sum is
\[
S_{H_0}(t)=h\,1_{\{h\mid t\}}.
\]
For $t,u\neq0$, the centered Fourier block and the adjacency axes
are
\[
T_{tu}=
\frac{S_{H_0}(t-u)+S_{H_0}(t)\overline{S_{H_0}(u)}}M\geq0,
\qquad
\beta_t=\frac{q-1}{M}S_{H_0}(t)\geq0.
\]
Since zero belongs to a smallest part, all triangle inequalities
in the lower-bound proof are equalities. Lift this labeling along
the quotient map $\mathbb Z_{aM}\to\mathbb Z_M$. The resulting
parts have sizes $a,\ldots,a,qa$ and the required nested coset
structure. Under this lift both Fourier norms are multiplied by
$a$, so their ratio is unchanged. The right side of
\eqref{eq:multipartite-sharp-bound} is likewise unchanged.
This proves sufficiency for every labeling described in the
statement.
\end{proof}

\subsection{Tripartite graphs and explicit families}

\begin{corollary}\label{cor:tripartite-sharp}
If $1\leq a\leq b\leq c$, then
\begin{equation}\label{eq:tripartite-sharp}
\operatorname{FR}_{\min}(K_{a,b,c})
\geq
\frac{4(b^2+bc+c^2)}
{(a+b+c)\sqrt{2(ab+ac+bc)}}.
\end{equation}
Equality holds precisely when $a=b=c$, or when $a=b$ and
$c=2ka$ for an integer $k\geq1$.
For every positive even integer $q$,
\begin{equation}\label{eq:tripartite-explicit}
\operatorname{FR}_{\min}(K_{a,a,qa})
=\frac{4(q^2+q+1)}{(q+2)\sqrt{4q+2}}.
\end{equation}
In particular,
\[
\operatorname{FR}_{\min}(K_{a,a,2a})=\frac7{\sqrt{10}},
\qquad
\operatorname{FR}_{\min}(K_{a,a,4a})=\frac{7\sqrt2}{3}.
\]
\end{corollary}

\begin{proof}
For three parts,
\[
N^2+(a^2+b^2+c^2)-2aN=2(b^2+bc+c^2).
\]
Apply Theorems~\ref{thm:multipartite-sharp-bound} and
\ref{thm:multipartite-sharp-equality}. Substitution of $b=a$ and
$c=qa$ gives \eqref{eq:tripartite-explicit}.
\end{proof}

\begin{remark}
In \eqref{eq:tripartite-explicit}, the two small parts project to
antipodal singleton labels in $\mathbb Z_{q+2}$.
The general equality family has an equally explicit value: if
$h\geq1$ and $q$ is a positive multiple of $h$, then
\begin{equation}\label{eq:multipartite-explicit}
\operatorname{FR}_{\min}
\left(K_{\underbrace{a,\ldots,a}_{h\text{ parts}},qa}\right)
=
\frac{2\left(h^2+2hq+2q^2-h-2q\right)}
{(h+q)\sqrt{h(h-1+2q)}}.
\end{equation}
For example,
$\operatorname{FR}_{\min}(K_{a,a,a,3a})=\sqrt6$.
Outside the equality families, the minimum is strictly larger
than \eqref{eq:multipartite-sharp-bound}; the theorem does not
otherwise determine its exact value.
\end{remark}

\subsection{Comparison with graph energy}

The energy lower bound can be calculated from a matrix with one
row for each part. This avoids finding all eigenvalues of the
$N$ by $N$ adjacency matrix. General background on graph spectra
can be found in \cite{BrouwerHaemers2012}.

For complete multipartite graphs, the identity
$\mathcal E(K_{n_1,\ldots,n_r})=2\rho$, where $\rho$ is the
spectral radius, and the equation \eqref{eq:rho} determining
$\rho$ are classical; see
\cite[Lemma~1.2 and the discussion following equation~(1.2)]
{StevanovicGutmanRehman2015}.
We recall their derivation below, then show that our lower bound
\eqref{eq:multipartite-sharp-bound} is strictly larger than the
resulting energy bound whenever the part sizes are not all equal.

\begin{proposition}\label{prop:multipartite-energy}
Let $\rho>0$ be the unique solution of
\begin{equation}\label{eq:rho}
\sum_{j=1}^r\frac{n_j}{\rho+n_j}=1.
\end{equation}
Then
\[
\mathcal E(K_{n_1,\ldots,n_r})=2\rho,
\qquad
\operatorname{FR}_{\min}(K_{n_1,\ldots,n_r})
\geq\frac{2\rho}{\sqrt{N^2-\sum_jn_j^2}}.
\]
If the part sizes are not all equal, the lower bound in
\eqref{eq:multipartite-sharp-bound} is strictly larger than this
energy bound.
\end{proposition}

\begin{proof}
Vectors whose coordinates sum to zero on each individual part
are annihilated by the adjacency matrix. Their orthogonal
complement is the space of vectors constant on each part.
In the orthonormal basis $1_{C_j}/\sqrt{n_j}$, the restriction
of the adjacency matrix to this space is
\[
Q=ww^*-\operatorname{diag}(n_1,\ldots,n_r),
\qquad w=(\sqrt{n_1},\ldots,\sqrt{n_r})^T.
\]
The quadratic form of $Q$ is negative definite on $w^\perp$,
so $Q$ has at most one positive eigenvalue. Moreover,
\[
w^*Qw=N^2-\sum_{j=1}^r n_j^2=2s>0,
\]
so it has exactly one positive eigenvalue, denoted by $\rho$.

For a corresponding nonzero eigenvector $z$,
\[
(\rho+n_j)z_j=\sqrt{n_j}\,w^*z.
\]
Here $w^*z\neq0$, since otherwise every $z_j$ would vanish.
Dividing by $\rho+n_j$, multiplying by $\sqrt{n_j}$, and
summing over $j$ gives
\[
w^*z=(w^*z)\sum_{j=1}^r\frac{n_j}{\rho+n_j}.
\]
Cancellation of $w^*z$ proves \eqref{eq:rho}. Its left side
is continuous and strictly decreasing for $\rho\geq0$, equals
$r>1$ at $\rho=0$, and tends to zero as $\rho\to\infty$.
Thus \eqref{eq:rho} has a unique positive solution.

Since $\operatorname{tr}Q=0$, the nonpositive eigenvalues
of $Q$ sum to $-\rho$, and their absolute values sum to $\rho$.
All remaining adjacency eigenvalues are zero. Hence the graph
energy is $2\rho$, and \eqref{eq:energy-bound} gives the
displayed lower bound.

The maximum vertex degree is $N-a$, so $\rho\leq N-a$.
If the part sizes are unequal, then
\[
\mu-a=\frac1N\sum_{j=1}^r n_j(n_j-a)>0.
\]
Consequently,
\[
N+\mu-2a=(N-a)+(\mu-a)>\rho.
\]
Dividing twice these quantities by $\sqrt{2s}$ proves that
the lower bound in \eqref{eq:multipartite-sharp-bound}
is strictly larger than the energy bound.

Strictness of the energy inequality itself also follows
from Corollary~\ref{cor:regular}: a vertex in part $j$
has degree $N-n_j$, so unequal part sizes make the graph
nonregular.
\end{proof}

\begin{corollary}\label{cor:multipartite-blowup}
For every positive integer $k$,
\[
\operatorname{FR}_{\min}(K_{kn_1,\ldots,kn_r})
\leq\operatorname{FR}_{\min}(K_{n_1,\ldots,n_r}).
\]
Equality holds whenever the part sizes satisfy
Theorem~\ref{thm:multipartite-sharp-equality}. In particular,
\[
\operatorname{FR}_{\min}
\left(K_{\underbrace{a,\ldots,a}_{r\text{ parts}}}\right)
=2\sqrt{1-\frac1r}.
\]
\end{corollary}

\begin{proof}
The first assertion follows from the periodic lift of
Theorem~\ref{thm:lift}. The sharp lower bound and its equality
conditions are unchanged by multiplying every part size by $k$.
This proves the second assertion. Substitution of $n_j=a$ in
\eqref{eq:multipartite-sharp-bound} gives the balanced value.
\end{proof}

The balanced value can also be understood directly through
blow-ups. The complete multipartite graph with $r$ parts of
size $a$ is $K_r^{(a)}$. The complete graph $K_r$ has energy
$2(r-1)$ and $r(r-1)/2$ edges, and its adjacency matrix
$J_r-I_r$ is diagonalized by the cyclic Fourier matrix.
It therefore attains the energy bound, with Fourier ratio
$2\sqrt{1-1/r}$. The periodic lift in Theorem~\ref{thm:lift}
gives a labeling of $K_r^{(a)}$ with the same Fourier ratio.
Since the energy bound is also unchanged under this blow-up,
this labeling is minimizing.

\subsection{A refinement from the degree sequence}

The zero diagonal of an adjacency matrix gives a further bound
which does not require the graph to be complete multipartite.

\begin{lemma}\label{lem:degree-fourier-bound}
Let $G$ be a simple undirected graph with $N$ vertices and $s>0$
edges. Write $\overline d=2s/N$ for its average degree and
\[
\delta=\max_{x\in V(G)}|\deg(x)-\overline d|.
\]
Then
\begin{equation}\label{eq:degree-fourier-bound}
\operatorname{FR}_{\min}(G)
\geq\frac{2\overline d+4\delta}{\sqrt{2s}}.
\end{equation}
\end{lemma}

\begin{proof}
For any labeling put $B=UAU^*$ and $\beta_t=B_{t0}$ for $t\neq0$.
Then $B_{00}=\overline d$. Because $A$ has zero diagonal, for each
$k\in\mathbb Z_N$,
\[
\sum_{u\in\mathbb Z_N}B_{u+k,u}
=\sum_{x\in\mathbb Z_N}A_{xx}e^{-2\pi ikx/N}=0.
\]
For $k=0$, the off-origin diagonal entries therefore have sum
$-\overline d$, so their absolute values sum to at least
$\overline d$. For $k\neq0$, the two terms on the Fourier axes
are $B_{k0}=\beta_k$ and $B_{0,-k}=\beta_k$; the latter identity
uses the reality and symmetry of $A$. Thus
\[
\sum_{\substack{u\neq0\\u+k\neq0}}B_{u+k,u}=-2\beta_k.
\]
The triangle inequality on each of these disjoint cyclic diagonals
gives
\[
\|B\|_{\ell^1}\geq
2\overline d+4\sum_{t\neq0}|\beta_t|.
\]
Finally, Fourier inversion gives
\[
\deg(x)-\overline d
=\sum_{t\neq0}\beta_t e^{2\pi itx/N},
\]
so $\sum_{t\neq0}|\beta_t|\geq\delta$.
Divide by $\sqrt{2s}$ and minimize over labelings.
\end{proof}

\begin{corollary}\label{cor:multipartite-degree-refinement}
With the multipartite notation above,
\begin{equation}\label{eq:multipartite-degree-refinement}
\operatorname{FR}_{\min}(K_{n_1,\ldots,n_r})
\geq
\frac{2(N-\mu)+4\max\{\mu-a,b-\mu\}}{\sqrt{2s}}.
\end{equation}
This coincides with \eqref{eq:multipartite-sharp-bound} when
$\mu-a\geq b-\mu$, and is stronger when $b-\mu>\mu-a$.
\end{corollary}

\begin{proof}
The average degree is $N-\mu$, while a vertex in part $j$ has
degree $N-n_j$. Thus the largest absolute degree deviation is
$\max\{\mu-a,b-\mu\}$. Apply Lemma~\ref{lem:degree-fourier-bound}.
\end{proof}

\begin{remark}
The equality classification in
Theorem~\ref{thm:multipartite-sharp-equality} concerns
\eqref{eq:multipartite-sharp-bound}. The refinement
\eqref{eq:multipartite-degree-refinement} is stated separately.
\end{remark}

For completeness, the graph identities
\[
K_{n_1,\ldots,n_{2q}}=
K_{n_1,n_2}\vee\cdots\vee K_{n_{2q-1},n_{2q}}
\]
and
\[
K_{n_1,\ldots,n_{2q+1}}=
K_{n_1,n_2}\vee\cdots\vee K_{n_{2q-1},n_{2q}}
\vee\overline K_{n_{2q+1}}
\]
express complete multipartite graphs as joins. Here $G\vee H$
adds every edge between disjoint copies of $G$ and $H$, and
$\overline K_m$ is an independent set of size $m$.
These identities do not by themselves give identities for
edge complexity. The Fourier formula
\eqref{eq:multipartite-formula} and the sharp bound
\eqref{eq:multipartite-sharp-bound} retain the information needed
for the labeling problem. Equation~\eqref{eq:rho} determines
the graph energy through Proposition~\ref{prop:multipartite-energy},
which also gives the corresponding energy lower bound.
General graph-energy estimates for joins are given in
\cite{GIIST2026}.

\section{Harmonic complexity of complete bipartite graphs}
\label{sec:harmonic}

We now apply the Fourier ratio to the eigenspaces of the Laplacian.
The orthogonal projector onto an eigenspace is a matrix indexed by
vertices, so we can ask the same labeling question as for the
adjacency matrix. For the combinatorial Laplacian
$L=\operatorname{diag}(\deg(v))-A$, let $\Pi_\lambda$ denote the
orthogonal projector onto the $\lambda$-eigenspace. Under a labeling
it is a function of two vertex variables. Following
\cite{GuptaIosevich2026}, define
$$
\operatorname{FR}_{\min}(\Pi_\lambda)
=\min_\pi\frac{\|\widehat{\Pi_{\lambda,\pi}}\|_1}{\|\widehat{\Pi_{\lambda,\pi}}\|_2}.
$$
Only nonzero projectors are considered. Their ranks are the
multiplicities of the corresponding eigenvalues.

The multiplicity lower bound for Laplacian projectors was proved
by Gupta and Iosevich \cite[Theorem~7.2]{GuptaIosevich2026}.
They also showed that the natural cyclic labeling of a circulant
graph attains this bound for every Laplacian projector
\cite[Theorem~7.3]{GuptaIosevich2026}.
The equality criterion for an arbitrary orthogonal projector
also appears in
\cite[Proposition~11.2]{GuptaEtAlAdditive2026}.
We recall the bound and give a short proof of this criterion
for a fixed cyclic labeling.

\begin{lemma}\label{lem:projector}
If $P$ is an orthogonal projector of rank $r>0$, then
$$
\operatorname{FR}(P)\geq\sqrt r.
$$
For a fixed cyclic labeling, equality holds if and only if $P$ is
circulant, that is, $P(x+1,y+1)=P(x,y)$ for all $x,y$.
\end{lemma}

\begin{proof}
The matrix $B=UPU^*$ is an orthogonal projector. Its diagonal entries
are nonnegative and sum to $\operatorname{tr}(B)=r$. Therefore
$$
\|\widehat P\|_1=\sum_{a,b}|B_{ab}|\geq r,
  \|\widehat P\|_2=\sqrt r.
$$
Equality holds exactly when every off-diagonal entry of $B$ vanishes.
If $B$ is diagonal, Fourier inversion gives
$$
P(x,y)=\frac{1}{N}\sum_{t=0}^{N-1}B_{tt}
e^{\frac{2\pi it(x-y)}{N}},
$$
which depends only on $x-y$. Conversely, if $P$ depends only on
$x-y$, summing in one variable shows that $UPU^*$ is diagonal.
\end{proof}

In particular, a projector with nonconstant diagonal cannot attain
the multiplicity bound under any labeling. This observation resolves
most of the equality questions for an unbalanced complete bipartite
graph before any Fourier sums need to be evaluated.

\begin{theorem}\label{thm:unbalanced-projectors}
Let $1\leq m<n$ and $N=m+n$. The nonzero Laplacian projectors of
$K_{m,n}$ satisfy
$$
\operatorname{FR}_{\min}(\Pi_0)=1,
\operatorname{FR}_{\min}(\Pi_N)=\frac{d_{m,N}^2}{mn}
\geq\frac{n}{m}>1.
$$
Equality in the lower bound
$\operatorname{FR}_{\min}(\Pi_N)\geq n/m$
holds exactly when $m\mid n$.
Since $\Pi_N$ has rank one, its multiplicity bound is $1$,
which is strictly smaller than $n/m$.
The labels minimizing $\operatorname{FR}(\Pi_N)$ are exactly those
minimizing the adjacency Fourier ratio. Moreover,
$$
\operatorname{FR}_{\min}(\Pi_m)>\sqrt{n-1},
$$
and, if $m>1$,
$$
\operatorname{FR}_{\min}(\Pi_n)>\sqrt{m-1}.
$$
Thus the zero-eigenvalue projector is the only one attaining its
multiplicity bound in the unbalanced case.
\end{theorem}

\begin{proof}
Let $C$ be the labels on the part of size $m$ and let
$B=\mathbb Z_N\setminus C$. Vectors supported on $C$ with coordinate
sum zero are eigenvectors of $L$ with eigenvalue $n$. This space has
dimension $m-1$. Vectors supported on $B$ with coordinate sum zero
have eigenvalue $m$ and form a space of dimension $n-1$. The
remaining two dimensions are the constant vector and the unit vector
\begin{equation}\label{eq:top-vector}
v=\sqrt{\frac{n}{mN}}1_C-\sqrt{\frac{m}{nN}}1_B
=\sqrt{\frac{N}{mn}}(1_C-\frac{m}{N}).
\end{equation}
Direct substitution into $L$ gives eigenvalues $0$ and $N$ for these
vectors, respectively. Their orthogonality and the dimension count
give the full spectral decomposition.

The constant projector is $\Pi_0(x,y)=\frac{1}{N}$, so its Fourier transform
has just one nonzero entry, equal to $1$. This proves the first
identity.

For the top projector $\Pi_N(x,y)=v(x)v(y)$, its two-dimensional
transform factors as $\widehat v(a)\widehat v(b)$. The vector in the equation
\eqref{eq:top-vector} has zero mean, and at nonzero frequencies
$$
\widehat v(t)=\frac{S_C(t)}{\sqrt{mn}}.
$$
Since $\|v\|_2=1$, the projector has Frobenius norm one. Hence
\begin{equation}\label{eq:top-projector-formula}
\operatorname{FR}(\Pi_{N,C})=
\|\widehat v\|_1^2=\frac{D_N(C)^2}{mn}.
\end{equation}
Taking the minimum, and applying Theorem~\ref{thm:coset}, gives the
formula, its lower bound, and the divisibility condition. Both the equation
\eqref{eq:top-projector-formula} and the adjacency formula increase
strictly with $D_N(C)$, proving the assertion about minimizing labels.

The projector onto the eigenvalue $m$ is
$$
\Pi_m(x,y)=1_B(x)1_{\{x\}}(y)-\frac{1}{n}1_B(x)1_B(y).
$$
Its diagonal is zero on $C$ and $\frac{n-1}{n}$ on $B$. This is
nonconstant because $n\geq2$ and both parts are nonempty. Similarly,
when $m>1$,
$$
\Pi_n(x,y)=1_C(x)1_{\{x\}}(y)-\frac{1}{m}1_C(x)1_C(y)
$$
has nonconstant diagonal. Lemma~\ref{lem:projector} excludes equality
for every labeling. Since the number of labelings is finite, the
minimum is strictly larger than the respective multiplicity bound.
\end{proof}

For balanced complete bipartite graphs, the projector values under
the parity labeling were given by Gupta and Iosevich
\cite[Corollary~7.4]{GuptaIosevich2026}. We recall this result,
together with a short proof, to contrast it with the unbalanced
case above.

\begin{theorem}\label{thm:balanced-projectors}
For $K_{m,m}$, the parity labeling simultaneously gives
$$
\operatorname{FR}_{\min}(\Pi_0)=1,
 \operatorname{FR}_{\min}(\Pi_{2m})=1,
$$
and, when $m\geq2$,
$$
\operatorname{FR}_{\min}(\Pi_m)=\sqrt{2m-2}.
$$
Every nonzero Laplacian projector therefore attains its multiplicity
bound under the same labeling.
\end{theorem}

\begin{proof}
Under the parity labeling the adjacency matrix depends only on the
difference of the labels and is circulant. The degree matrix is
$mI$, so the Laplacian is also diagonal in the cyclic Fourier basis.
Its spectral projectors are therefore circulant. The constant
eigenspace and the eigenspace at $2m$ have rank one. The remaining
eigenvalue $m$ has rank $2m-2$, and is absent when $m=1$. Apply
Lemma~\ref{lem:projector}.
\end{proof}

\begin{remark}
When $m=n$, the two spaces supported on the separate parts have the
same Laplacian eigenvalue. Their projectors must be added to obtain
the full spectral projector. Thus the two intermediate projectors
of the unbalanced case become a single projector when $m=n$.
\end{remark}

For the intermediate projectors in the unbalanced case, their exact
minima remain to be determined. The formulas above give a direct
starting point. For example, when $m>1$,
$$
\widehat{\Pi_{n,C}}(a,b)=\frac{1}{N}
(S_C(a+b)-\frac{S_C(a)S_C(b)}{m}).
$$
The formula for $\Pi_m$ is obtained by replacing $C,m$ by $B,n$.
Unlike the top-projector formula, these expressions retain information
about the phases and the additive arrangement of the frequencies.

\section{Recovering connections between stores and warehouses}
\label{sec:recovery}

Suppose that a company has $n$ stores and $m$ warehouses. A store
can be supplied by a warehouse when the drive between them takes at
most one day. We represent the possible deliveries by a bipartite
graph. Its two parts are the stores and the warehouses, and an edge
records an allowed delivery. The question is whether we must check
every one of the $mn$ store--warehouse pairs to obtain an accurate
picture of the graph.

A warehouse may be too far from some stores, so this graph need
not be complete. The formula in
Theorem~\ref{thm:sharp} applies to $K_{m,n}$, for which every store
can be supplied by every warehouse. For an incomplete graph we need
a bound for the Fourier ratio of its own adjacency function.
Divisibility of $n$ by $m$ does not supply such a bound.

Let $N=m+n$, and fix a labeling by $\mathbb Z_N$. Let $W$ be the
warehouse labels and let $S$ be the store labels. Write $s>0$ for
the number of edges. As elsewhere in the paper, $f(x,y)=1$ when
$x$ and $y$ are adjacent, and $f(x,y)=0$ otherwise. We keep the
full symmetric adjacency function on $\mathbb Z_N^2$. Each edge
appears twice, and $f(x,y)^2=f(x,y)$. Therefore
\begin{equation}\label{eq:recovery-norm}
\|f\|_2^2=\sum_{x,y\in\mathbb Z_N}f(x,y)^2=2s.
\end{equation}
Thus the norm is known exactly whenever the edge count is known.
For $K_{m,n}$ it is $\sqrt{2mn}$. For an incomplete graph, $m$ and
$n$ alone do not determine $s$.

\subsection{From sampled values to an approximation}

Bounds for the number of Fourier frequencies needed to
approximate an adjacency function are given in
\cite[Theorem~5.8]{GuptaEtAlAdditive2026}.
These concern representation size. Recovery from sampled
entries requires a separate observation argument; see
\cite[Remark~5.9]{GuptaEtAlAdditive2026}.
Choose a set of $q$ distinct ordered pairs uniformly from all
$q$-element subsets of $\mathbb Z_N^2$. List them as
$z_1,\ldots,z_q$ and observe $f(z_j)$. A pair of stores or a pair
of warehouses has value zero without a driving-time calculation.
The two orientations of an edge also share one calculation. Thus
there are at most $q$ actual driving-time checks. We sample from
the whole square $\mathbb Z_N^2$, including its known zero entries.

Among the functions agreeing with these observations, choose one
whose Fourier coefficients have the smallest sum of absolute
values. In formulas, choose a minimizer $g$ of
\begin{equation}\label{eq:recovery-minimization}
\min_g\|\widehat g\|_1
 \text{subject to }g(z_j)=f(z_j),  1\leq j\leq q.
\end{equation}

This is the basis-pursuit optimization principle
\cite{ChenDonohoSaunders1998}, applied to the Fourier coefficients.
The constraints prescribe the sampled values of their inverse
Fourier transform.

The unknown values of $g$ are allowed to be complex numbers at this
stage. We will turn them into zeros and ones after estimating the
error. The true function $f$ is one of the allowed choices, so
every minimizer satisfies
$$
\|\widehat g\|_1\leq\|\widehat f\|_1.
$$
This inequality will control the Fourier transform of the error.

The Fourier-ratio recovery theorem of Burstein, Iosevich, and Nathan
\cite{BursteinIosevichNathan2026} explains how this minimization
recovers a function with small Fourier ratio. Burstein, Iosevich,
and Krause \cite{BursteinIosevichKrause2026} improve the sampling
estimate used in that argument. We use their result to obtain the
following bound. The proof also explains why the sampling theorem
must hold for a whole class of functions at once.

\begin{theorem}\label{thm:graph-approximation}
Let $f$ be the adjacency function of a graph on $\mathbb Z_N$ with
$s>0$ edges, and suppose that $\operatorname{FR}(f)\leq R$, where
$R\geq1$. Fix $0<\varepsilon<1$ and put
$$
k=\lceil\frac{4R^2}{\varepsilon^2}\rceil.
$$
There is an absolute constant $C>0$ such that, if $k\leq N^2$ and
\begin{equation}\label{eq:recovery-samples}
Ck\log(\frac{2eN^2}{k})\log(2k)
\leq q\leq\frac{N^2}{2},
\end{equation}
then a uniformly chosen set of $q$ distinct ordered pairs has the
following property with probability at least $1-CN^{-2}$.
Every minimizer in the equation~\eqref{eq:recovery-minimization}
satisfies
$$
\|g-f\|_2\leq\varepsilon\sqrt{2s}.
$$
The labeling and the bound $R$ are fixed before the sample is drawn.
\end{theorem}

The proof uses the complex-matrix sampling theorem of
Burstein, Iosevich, and Krause
\cite[Theorem~1.6 and Corollary~5.2]{BursteinIosevichKrause2026}.
For clarity, the substitutions into their result are listed below.
Here $h$ denotes a function on $\mathbb Z_N^2$, and
$\Omega=\{z_1,\ldots,z_q\}$ is the set of sampled ordered pairs.

\begin{center}
\begin{tabular}{@{}lll@{}}
\toprule
Quantity & Notation in the source & Present application \\
\midrule
Ambient dimension & $n$ & $N^2$ \\
Number of sampled rows & $m$ & $q$ \\
Norm-class parameter & $r$ & $k=\lceil4R^2/\varepsilon^2\rceil$ \\
Coefficient vector & $y$ & $\widehat h$ \\
Sampled row set & $S$ & $\Omega$ \\
Relative energy distortion & $\varepsilon_{\mathrm{src}}$ & $1/8$ \\
Failure exponent & $A_0$ & $1$ \\
\bottomrule
\end{tabular}
\end{center}

The source writes its energy distortion as $\varepsilon$; we have
added the subscript in the table to distinguish it from the
recovery accuracy in Theorem~\ref{thm:graph-approximation}.
The distortion is fixed at $1/8$, while the recovery accuracy
determines $k$. Likewise, the source's $m$ and $n$ denote a sample
count and an ambient dimension, not the two part sizes of our graph.
Corollary~5.2 of the cited paper gives the uniform sampling of
exactly $q$ distinct rows required here.

\begin{proof}
We apply the complex-matrix sampling theorem of
\cite{BursteinIosevichKrause2026} to the inverse Fourier transform
on $\mathbb Z_N^2$. To check the normalization, consider the matrix
$$
H_{z,\xi}=e^{\frac{2\pi iz\cdot\xi}{N}},
z,\xi\in\mathbb Z_N^2.
$$
Its entries have absolute value one, and orthogonality of the
characters gives $H^*H=N^2I$. Fourier inversion reads
$h=\frac{1}{N}H\widehat h$. Thus the matrix theorem applies with
dimension $N^2$.
To make the passage from vectors to functions explicit, write
$H_\Omega$ for the restriction of $H$ to the sampled rows.
The lower energy estimate in the cited result has the form
\[
(1-\varepsilon_{\mathrm{src}})\|y\|_2^2
\leq\frac1q\|H_\Omega y\|_2^2,
\]
simultaneously for all coefficient vectors $y$ satisfying
$\|y\|_1\leq\sqrt{k}\|y\|_2$.
Here $\varepsilon_{\mathrm{src}}$ denotes the energy distortion,
as in the table above.

We take $y=\widehat h$. Parseval's identity gives
$\|y\|_2=\|h\|_2$, so the condition on the coefficient vector
becomes $\|\widehat h\|_1\leq\sqrt{k}\|h\|_2$.
Moreover, Fourier inversion gives $H\widehat h=Nh$, and hence
\[
\frac1q\|H_\Omega\widehat h\|_2^2
=\frac1q\sum_{j=1}^q|(H\widehat h)(z_j)|^2
=\frac{N^2}{q}\sum_{j=1}^q|h(z_j)|^2.
\]
This explains the normalization of the sampled energy.
The remaining hypotheses of the sampling result hold as well.
We have $1\leq k\leq N^2$, and under the substitutions
$n=N^2$, $m=q$, and $r=k$, its sample-size condition becomes
\eqref{eq:recovery-samples}, provided $C$ is chosen sufficiently
large. For the choices $\varepsilon_{\mathrm{src}}=1/8$ and
$A_0=1$ in the table, the source's constant
$C_{\varepsilon_{\mathrm{src}},A_0}$ is an absolute constant.
In particular, it is independent of our recovery accuracy
$\varepsilon$, which enters the sample bound through $k$.
The source's failure bound $Cn^{-A_0}$ becomes $CN^{-2}$.
Taking its energy error to be $\frac{1}{8}$ and
its failure exponent to be $1$, we obtain, with probability at
least $1-CN^{-2}$,
\begin{equation}\label{eq:recovery-sampling-estimate}
\frac{7}{8}\|h\|_2^2
\leq\frac{N^2}{q}\sum_{j=1}^q|h(z_j)|^2
\end{equation}
simultaneously for every $h$ satisfying
$\|\widehat h\|_1\leq\sqrt k\|h\|_2$. The factor $\frac{N^2}{q}$
puts the sampled sum on the scale of the sum over all pairs.

Now let $h=g-f$. It vanishes at every sampled pair. Moreover,
the minimizing property of $g$ and the triangle inequality imply
$$
\|\widehat h\|_1
\leq\|\widehat g\|_1+\|\widehat f\|_1
\leq2\|\widehat f\|_1
\leq2R\|f\|_2.
$$
Suppose that $\|h\|_2>\varepsilon\|f\|_2$. Then
$$
\|\widehat h\|_1
<\frac{2R}{\varepsilon}\|h\|_2
\leq\sqrt k\|h\|_2,
$$
so the sampling estimate applies to $h$. Its right side is zero,
forcing $h=0$, a contradiction. Hence
$\|h\|_2\leq\varepsilon\|f\|_2$. The equation
\eqref{eq:recovery-norm} completes the proof.
\end{proof}

Although $g$ depends on the sampled values, the event used in the
proof controls every $h$ in the stated class. This is why we can
apply the estimate to the unknown error $g-f$. The factor
$\frac{7}{8}$ remains fixed as we ask for greater accuracy.
The accuracy enters through $k$, which is proportional to
$\frac{R^2}{\varepsilon^2}$.

For comparison, the estimate of Haviv and Regev
\cite{HavivRegev2015} used in the earlier recovery argument leads,
up to absolute constants, to the sample count
$k\log^2(2k)\log(N^2)$. The equation~\eqref{eq:recovery-samples}
removes one factor of $\log(2k)$ and replaces $\log(N^2)$ by
$\log(\frac{2eN^2}{k})$. If $k>N^2$, or if the lower bound for
$q$ exceeds $\frac{N^2}{2}$, we may instead check all $mn$
store--warehouse pairs and recover the graph exactly.

The minimization uses only the observed values. It does not require
the unknown edge count $s$. If $s$ is supplied separately, we may
also require $\sum_{x,y}g(x,y)=2s$. We may likewise require that
$g$ be real and symmetric, that $0\leq g(x,y)\leq1$, and that it
vanish on pairs within either part. These conditions all hold for
$f$. The same proof therefore applies with these additional
conditions, because it uses only agreement at the sampled points
and the inequality $\|\widehat g\|_1\leq\|\widehat f\|_1$.

\subsection{Counting the incorrectly recovered connections}

An approximation by arbitrary numbers does not yet give a graph.
The fact that the true values are zero or one provides a simple
way to finish. First set
$$
u(x,y)=\frac{\operatorname{Re}g(x,y)+\operatorname{Re}g(y,x)}{2}
$$
when one label lies in $W$ and the other in $S$, and set $u(x,y)=0$
otherwise. Taking real parts, averaging the two orientations, and
setting the known zero entries to zero cannot increase the squared
error. For the averaging step this follows from
$$
|\frac{a+b}{2}-t|^2
\leq\frac{|a-t|^2+|b-t|^2}{2},
$$
applied with $t=f(x,y)=f(y,x)$.
More explicitly, for each $x\in W$ and $y\in S$, symmetry of $u$
and the preceding inequality give
\[
\begin{aligned}
|u(x,y)-t|^2+|u(y,x)-t|^2
&=2|u(x,y)-t|^2\\
&\leq|\operatorname{Re}g(x,y)-t|^2
     +|\operatorname{Re}g(y,x)-t|^2\\
&\leq|g(x,y)-t|^2+|g(y,x)-t|^2.
\end{aligned}
\]
The last inequality uses the fact that $t$ is real: discarding
the imaginary part cannot increase the squared distance to $t$.
Summing over $W\times S$ includes both orientations of every
possible store--warehouse connection. On all remaining ordered
pairs, both $u$ and $f$ are zero.
It follows that
$$
\|u-f\|_2\leq\|g-f\|_2.
$$

Join a store and a warehouse in the recovered graph when
$u(x,y)\geq\frac{1}{2}$. Denote its edge set by $E_*$. We count
an error whenever an edge has been omitted or a nonexistent edge
has been inserted. The set of such errors is the symmetric
difference $E\triangle E_*$.

If an edge is omitted, its true value is $1$ and its approximating
value is less than $\frac{1}{2}$. If a nonexistent edge is inserted,
its true value is $0$ and its approximating value is at least
$\frac{1}{2}$. In either case, the absolute error is at least
$\frac{1}{2}$ in each orientation. Thus every incorrect connection
contributes at least $\frac{1}{2}$ to $\|u-f\|_2^2$.
The contribution is $2(1/2)^2=1/2$ because there are two
orientations, each with squared error at least $1/4$.
The threshold case causes no exception: if $u(x,y)=1/2$, the
recovered graph includes the edge, so an error occurs only
when the true value is zero, and then the squared error is
exactly $1/4$ in each orientation.

Summing over the incorrectly recovered connections and using
the approximation estimate gives
\[
\frac12|E\triangle E_*|
\leq\|u-f\|_2^2
\leq\|g-f\|_2^2
\leq\varepsilon^2\|f\|_2^2
=2\varepsilon^2s.
\]
Here the last equality uses \eqref{eq:recovery-norm}.
We obtain
\begin{equation}\label{eq:recovery-edge-errors}
|E\triangle E_*|
\leq2\|u-f\|_2^2
\leq4\varepsilon^2s.
\end{equation}
For example, relative $L^2$ error at most $\frac{1}{20}$ gives at
most $\frac{s}{100}$ incorrect connections. This is a count of all
errors, including both inserted and omitted edges.

\begin{corollary}\label{cor:graph-recovery}
Suppose that a bipartite graph with $m$ warehouses, $n$ stores,
and $s>0$ edges has a known labeling for which
$\operatorname{FR}(f)\leq R$, where $R\geq1$. Let $N=m+n$,
fix $0<\delta<1$, and put
$$
k_\delta=\lceil\frac{16R^2}{\delta}\rceil.
$$
There is an absolute constant $C>0$ such that, if
$k_\delta\leq N^2$ and
\begin{equation}\label{eq:recovery-hamming-samples}
Ck_\delta\log(\frac{2eN^2}{k_\delta})\log(2k_\delta)
\leq q\leq\frac{N^2}{2},
\end{equation}
a uniformly chosen set of $q$ distinct ordered pairs suffices to
construct an edge set $E_*$ satisfying
$$
|E\triangle E_*|\leq\delta s
$$
with probability at least $1-CN^{-2}$.
\end{corollary}

\begin{proof}
Take $\varepsilon=\frac{\sqrt\delta}{2}$ in
Theorem~\ref{thm:graph-approximation}. Its parameter $k$ becomes
$k_\delta$, and the equation~\eqref{eq:recovery-edge-errors}
bounds the number of incorrect edges by $\delta s$.
Indeed, $0<\delta<1$ implies
$0<\varepsilon=\sqrt{\delta}/2<1/2$, so this choice satisfies
the accuracy hypothesis of Theorem~\ref{thm:graph-approximation}.
The substitutions are
\[
k=\left\lceil\frac{4R^2}{\varepsilon^2}\right\rceil
 =\left\lceil\frac{4R^2}{\delta/4}\right\rceil
 =\left\lceil\frac{16R^2}{\delta}\right\rceil
 =k_\delta,
\qquad
4\varepsilon^2s=4\frac{\delta}{4}s=\delta s.
\]
Thus \eqref{eq:recovery-hamming-samples} is exactly the
sample-size condition required by that theorem.
Symmetrization and rounding are deterministic operations:
whenever the approximation estimate holds, the resulting edge
set satisfies the stated bound. Consequently, no additional
failure probability is introduced, and the success probability
remains at least $1-CN^{-2}$.
\end{proof}

Since $s\leq mn$, this also bounds the number of incorrect
connections by $\delta mn$. The estimate concerns the total number
of errors. It does not guarantee the same accuracy separately for
each warehouse or each store. The number of driving-time checks
can be small even though writing down the entire recovered graph
still requires $mn$ entries.

\subsection{What the complete bipartite calculation tells us}

Assume $m\leq n$ and $m\mid n$, and put
$$
R_{m,n}=\frac{4n^2}{(m+n)\sqrt{2mn}}.
$$
Theorem~\ref{thm:sharp} gives this ratio for a complete bipartite
graph when the warehouse labels form a coset of the subgroup of
order $m$. The quantity entering the sample bound is its square,
$$
R_{m,n}^2=\frac{8n^3}{m(m+n)^2}.
$$
Writing $r=\frac{n}{m}$, we have
$$
R_{m,n}^2=\frac{8r^3}{(1+r)^2}.
$$
Since $r\leq1+r\leq2r$ for $r\geq1$, it follows that
$$
2r\leq R_{m,n}^2\leq8r.
$$

Suppose now that the actual delivery graph, possibly incomplete,
has a known labeling satisfying
$$
\operatorname{FR}(f)\leq A R_{m,n}
$$
for a fixed constant $A\geq1$. To apply
Corollary~\ref{cor:graph-recovery}, put
$$
k_A=\lceil\frac{16A^2R_{m,n}^2}{\delta}\rceil.
$$
The bounds for $R_{m,n}^2$ show that
$$
\frac{32A^2n}{m\delta}\leq k_A
\leq\frac{129A^2n}{m\delta}.
$$
Here the upper estimate also accounts for the ceiling, since
$\frac{A^2n}{m\delta}\geq1$.

Suppose $k_A\leq N^2$. Its upper bound controls $k_A\log(2k_A)$,
while its lower bound controls $\log(\frac{2eN^2}{k_A})$.
Substitution into the corollary therefore shows that, for a
sufficiently large constant $C_A$ depending only on $A$, the condition
\begin{equation}\label{eq:warehouse-sample-count}
C_A\frac{n}{m\delta}
\log(2+\frac{n}{m\delta})
\log(2+\frac{m\delta(m+n)^2}{n})
\leq q\leq\frac{(m+n)^2}{2}
\end{equation}
suffices for at most $\delta s$ incorrect connections, with
probability at least $1-C(m+n)^{-2}$. If the sampling conditions
cannot be met, checking all $mn$ pairs still gives exact recovery.

For fixed $\frac{n}{m}$ and fixed $\delta$, the sufficient number
of random samples grows only logarithmically with $m+n$. The
estimate is useful when this count is much smaller than $mn$.
The constants are not explicit, so the formula describes how the
number of samples grows with the parameters.

For a complete bipartite graph with known parts, all connections
are already determined. Its Fourier ratio supplies a value against
which to compare an incomplete delivery graph. To use the bound in the equation
\eqref{eq:warehouse-sample-count} for an incomplete network, one
must establish the displayed bound for its actual adjacency
function. The optimal value over all labelings also does not by
itself provide a labeling that can be found before the graph is
known.

Finally, small Fourier ratio and binary values do not ensure exact
recovery from few point evaluations. A simple example explains
the obstruction. Start with $K_{t,t}$, where $t\geq2$, under the
parity labeling, and delete one unknown edge. Let $f_0$ and $f$ be
the original and modified adjacency functions. Since
$\operatorname{FR}(f_0)=\sqrt2$ and $\|f_0\|_2=\sqrt2t$, we have
$\|\widehat f_0\|_1=2t$. The Fourier transform of a function equal
to one at a single ordered pair and zero elsewhere has $\ell^1$
norm $2t$. Deleting an undirected edge changes two entries, so
the triangle inequality gives
$$
\operatorname{FR}(f)
\leq\frac{2t+4t}{\sqrt{2(t^2-1)}}
\leq2\sqrt6<5.
$$
Thus the Fourier ratio stays bounded, and the norm
$\sqrt{2(t^2-1)}$ is known exactly.

There are nevertheless $t^2$ possible locations for the deleted
edge. Choose that location uniformly at random. If at most $q<t^2$
distinct store--warehouse pairs are checked, the probability of
finding it among those pairs is at most $\frac{q}{t^2}$. If it is
not found, a single guess for its location contributes at most
$\frac{1}{t^2}$ to the overall success probability. Hence
the probability of exact recovery is at most
$\frac{q+1}{t^2}$. The same argument allows the checked pairs to
depend on the answers already obtained. A success probability
bounded away from zero therefore requires a number of checks
proportional to $t^2$. An approximation allowed one
incorrect connection can simply return $K_{t,t}$.

The distinction between recovering most connections and recovering
every connection is therefore necessary. To obtain an unconditional
sample bound for a delivery network, it remains to estimate the
Fourier ratio from the geometry of the store and warehouse locations.

\section{Further questions}
\label{sec:questions}

The complete classification for two- and three-element partite
sets leaves several distinct problems. Some concern the additive
structure of larger minimizing sets. Others concern the effect of
changing the ambient cyclic group, the interaction between several
parts, and the use of a minimizing labeling in recovery.

\subsection{Larger partite sets and additive energies}

Theorem~\ref{thm:coset} determines the minimum whenever $m\mid N$.
For three elements the remaining minimizers are precisely the
progressions with invertible common difference. The example
$K_{7,9}$ shows that this description cannot extend to all larger
partite sets, even when $m$ and $N$ are relatively prime.

\begin{question}
Determine $d_{m,N}$ and all its minimizing sets when
$4\leq m\leq N/2$ and $m\nmid N$. In particular, what are the
minimizing four-element sets? For a fixed $m\geq4$, does the interval
$I_m$ minimize $D_p$ for every sufficiently large prime $p$?
When it does, are all minimizers invertible affine images of $I_m$?
\end{question}

For composite $N$, unions of cosets and small modifications of
subgroups provide additional candidates. A classification should
explain how these compete with progressions. Even the next
cardinality has distinct arithmetic cases: when $N=2q$ with
$q$ odd, inverse images of minimizing two-element subsets of
$\mathbb Z_q$ provide natural four-element competitors in
$\mathbb Z_N$.

The three-element argument also raises a methodological question.
Its decisive step is a restriction on the fourth additive energy
when the third energy is minimal. For larger $m$, can polynomial
minorants of degree bounded in terms of $m$, together with
constraints on jointly possible additive energies, reduce the
problem to a small collection of candidate configurations?
The surviving configurations would still require separate
arithmetic comparisons. Indeed, a fixed number of moments cannot
distinguish progressions whose common differences have different
sufficiently large additive orders, even though their Fourier
norms can differ.

Inspired in this setting by the descent method of Burstein,
Iosevich, and Krause \cite{BursteinIosevichKrause2026}, one may
try to find a minimizing partite set by making successive small
changes that decrease its Fourier objective. A \emph{single
exchange} from $A\subset\mathbb Z_N$ removes one element of $A$
and inserts one element of $\mathbb Z_N\setminus A$. Write
\[
\mathcal N(A)=
\left\{(A\setminus\{x\})\cup\{y\}:
x\in A,\quad y\in\mathbb Z_N\setminus A\right\}
\]
for the resulting sets. The following example shows that descent
by single exchanges can stop at a set which is not globally
minimizing, even when $|A|=4$.

\begin{proposition}[Computer-assisted]
\label{prop:single-exchange-obstruction}
For four-element subsets $A\subset\mathbb Z_{17}$, put
$\Phi(A)=D_{17}(A)/17$, and let
\[
C=\{0,1,4,5\},\qquad P=\{0,1,2,3\}.
\]
Then
\[
\Phi(P)<1.3223,\qquad
1.3617<\Phi(C)<1.3618,
\]
whereas every $C'\in\mathcal N(C)$ satisfies
\[
\Phi(C')>1.3947.
\]
Thus $C$ is a strict local minimum under single exchanges but
is not a global minimum. Replacing $\{4,5\}$ by $\{2,3\}$
gives a decreasing two-element exchange.
\end{proposition}

\begin{proof}
There are exactly $4(17-4)=52$ distinct sets in $\mathcal N(C)$.
The accompanying verifier evaluates $P$, $C$, and all these
neighbors by exact rational arithmetic and integer square roots.
Appendix~\ref{app:single-exchange-certificate}
justifies the error bounds. The resulting certified enclosures,
rounded outward to nine decimal places, are
\begin{center}
\begin{tabular}{@{}lll@{}}
\toprule
Quantity & Lower endpoint & Upper endpoint\\
\midrule
$\Phi(P)$ & $1.322243933$ & $1.322243957$\\
$\Phi(C)$ & $1.361709120$ & $1.361709145$\\
$\displaystyle\min_{C'\in\mathcal N(C)}\Phi(C')$
  & $1.394738020$ & $1.394738045$\\
\bottomrule
\end{tabular}
\end{center}
The upper endpoint for $\Phi(P)$ is below the lower endpoint
for $\Phi(C)$, and the upper endpoint for $\Phi(C)$ is below
the lower endpoint for the minimum over its neighbors. The
enclosures also imply the stated decimal thresholds. The program
checks the underlying rational inequalities before rounding the
endpoints for display. Finally,
$(C\setminus\{4,5\})\cup\{2,3\}=P$, which proves the assertion
about two-element exchanges.
\end{proof}

By Theorem~\ref{thm:reduction}, the Fourier ratio of the labeling
of $K_{4,13}$ determined by $A$ is a strictly increasing function
of $\Phi(A)$. Hence $C$ also gives a strict local minimum of the
graph Fourier ratio under single exchanges of partite labels,
although $P$ gives a smaller ratio. No assertion that $P$ is a
global minimizer is needed for this conclusion. The obstruction
concerns the single-exchange strategy considered here; it does
not contradict the entropic descent results of
\cite{BursteinIosevichKrause2026}. Thus a general method for
finding a minimizing labeling must allow other moves, permit
temporary increases in the objective, or impose hypotheses that
exclude such local minima.

\subsection{Separation from the minimizing triples}

The classification of minimizers does not determine the next
smallest value. For $N\geq7$ with $3\nmid N$, let
\[
\mathcal M_N=
\left\{
t+u\{0,1,2\}:
t\in\mathbb Z_N,\quad u\in\mathbb Z_N^\times
\right\},
\]
and define
\[
\Delta_N=
\min_{\substack{C\subset\mathbb Z_N,\ |C|=3\\
                C\notin\mathcal M_N}}
\left(D_N(C)-D_N(\{0,1,2\})\right).
\]

\begin{question}
What is the exact value of $\Delta_N$, and which sets attain it?
How do the next minimizing affine classes depend on the
factorization of $N$?
\end{question}

The proof for triples already separates nonprogressions from
$\{0,1,2\}$ by at least $cN$ for an absolute $c>0$ and all
sufficiently large $N$ with $3\nmid N$. At prime orders this
applies to every nonminimizer. At composite orders, however,
progressions with noninvertible common difference can be much
closer to the minimum. For example, when $N=2k$ and $3\nmid N$,
the estimates for $W$ give, as $k\to\infty$,
\[
0<D_N(\{0,2,4\})-D_N(\{0,1,2\})
=N\bigl(W(k)-W(N)\bigr)=O(N^{-1}).
\]
Thus the arithmetic of $N$ is essential to a quantitative
description of near-minimizing labelings.

\subsection{Repeated blow-ups}

For a graph $G$ with at least one edge, put
\[
\beta(G)=
\inf_{k\geq1}\operatorname{FR}_{\min}(G^{(k)}).
\]
Periodic lifting makes the values nonincreasing when the blow-up
factor is replaced by a multiple. Consequently,
\[
\lim_{j\to\infty}
\operatorname{FR}_{\min}(G^{(j!)})=\beta(G).
\]
The energy bound gives
\[
\frac{\mathcal E(G)}{\sqrt{2|E(G)|}}
\leq\beta(G)\leq\operatorname{FR}_{\min}(G).
\]

\begin{question}
Is $\beta(G)$ always attained by a finite blow-up? If $G$ does
not attain its graph-energy bound, must $\beta(G)$ remain strictly
above that bound? Which graphs satisfy
$\beta(G)=\operatorname{FR}_{\min}(G)$?
\end{question}

Theorem~\ref{thm:blowup-converse} excludes equality with the energy
bound for each finite blow-up, but does not by itself give a gap
uniform in the blow-up factor. For unbalanced complete multipartite
graphs, Theorem~\ref{thm:multipartite-sharp-bound} does give such
a gap. The equality families in
Theorem~\ref{thm:multipartite-sharp-equality} also show that
attainment of the energy bound is not necessary for every
blow-up to preserve edge complexity.

For bipartite graphs, divisibility of $n$ by $m$ is sufficient
for this preservation. The more detailed question at a specified
blow-up factor remains:

\begin{question}
For which $1\leq m\leq n$ and $k\geq2$ does
\[
\operatorname{FR}_{\min}(K_{km,kn})
=\operatorname{FR}_{\min}(K_{m,n})
\]
hold? Which minimizing partite sets are inverse images under a
proper cyclic quotient map? In particular, is $m\mid n$ necessary
for the displayed equality to hold for every $k$?
\end{question}

The quotient question asks both when a periodic lift attains the
minimum and which minimizing sets admit such a representation.
The strict-decrease examples show that dividing both part sizes
by a common divisor does not determine the minimum in general.

\subsection{Multipartite minima and the degree bound}

The equality classification for the sharp multipartite bound
determines several infinite families exactly. Outside those
families, the Fourier sums of the parts interact through their
phases in \eqref{eq:multipartite-formula}. Minimizing the separate
quantities $D_j$ does not resolve this partition problem.
The upper bound \eqref{eq:multipartite-upper} provides a way to
compare candidate partitions with the lower bounds, but sharper
control of the cancellation between parts is needed.

\begin{question}
Determine the minimum and the minimizing partitions for further
tripartite families. In particular, what is
$\operatorname{FR}_{\min}(K_{a,a,qa})$ when $q\geq3$ is odd?
Is a periodic lift of a minimizing labeling of $K_{1,1,q}$
always sufficient to attain the minimum, or can a nonperiodic
labeling improve on every such lift?
\end{question}

There is a separate equality problem for the degree estimate
\[
\operatorname{FR}_{\min}(G)
\geq\frac{2\overline d+4\delta}{\sqrt{2s}},
\qquad
\delta=\max_x|\deg(x)-\overline d|.
\]
Which graphs and labelings attain this estimate? For complete
multipartite graphs, particular interest lies in the regime
$b-\mu>\mu-a$, where this estimate strictly improves the bound
whose equality cases have been classified. Can more information
from the degree sequence produce further sharp improvements?

\subsection{Laplacian projectors and common minimizing labelings}

For an unbalanced complete bipartite graph, the top projector and
the adjacency matrix have exactly the same minimizing labelings.
The intermediate projectors involve Fourier phases as well as
absolute values, so their optimization need not follow the same
principle.

\begin{question}
Determine $\operatorname{FR}_{\min}(\Pi_m)$ and, when $m>1$,
$\operatorname{FR}_{\min}(\Pi_n)$ for $K_{m,n}$ with $m<n$,
together with all minimizing labelings. For each intermediate
projector, when do its minimizing labelings also minimize the
adjacency Fourier ratio? For which part sizes can one labeling
simultaneously minimize adjacency and every nonzero Laplacian
projector ratio?
\end{question}

The coincidence with adjacency already fails for one intermediate
projector of $K_{2,7}$. Indeed, for
$C=\{0,d\}\subset\mathbb Z_N$ and $n=N-2>2$, the projector
$\Pi_n=vv^*$, with
$v=(1_{\{0\}}-1_{\{d\}})/\sqrt2$, satisfies
\[
\operatorname{FR}(\Pi_n)
=\frac{2N}{q^2}\cot^2\!\left(\frac{\pi}{2q}\right),
\qquad q=\frac{N}{\gcd(d,N)}.
\]
This follows by summing $|\sin(\pi td/N)|$, and the expression
is strictly increasing in $q$. For $N=9$, this projector is
therefore minimized at $q=3$, whereas adjacency is minimized
at $q=9$.

More generally, when individual spectral projectors have different
minimizing labelings, what tradeoffs are possible among their
Fourier ratios under a single labeling? The circulant equality
criterion settles simultaneous attainment of all multiplicity
bounds under a fixed labeling, but does not settle simultaneous
minimization when some of those bounds are strict.

\subsection{The choice of labeling group}

The cyclic group provides more than a convenient ordering: its
subgroups enter both the equality cases and the blow-up examples.
One can instead label the vertices by another finite abelian group
$H$ of the same order and use the unitary character transform on
$H\times H$. Even the two-element calculation depends on the
orders available for differences of labels; these orders need not
be the same in two groups of equal cardinality.

\begin{question}
How can the minimum Fourier ratio and its minimizing labelings
be described in terms of the structure of $H$? In particular,
which abelian groups of order $m+n$ minimize the resulting
edge complexity of $K_{m,n}$? Which analogues of the small-part
classifications and the blow-up descent theorem survive this
change of labeling group?
\end{question}

\subsection{Finding labelings and recovering graphs}

The recovery theorem assumes a known labeling with a controlled
Fourier ratio. A small minimum over all labelings does not supply
such a labeling before the graph is observed. Moreover, deleting
one unknown edge shows that a bounded ratio alone cannot ensure
exact recovery from a small number of adjacency queries.

\begin{question}
Which structural assumptions on an incomplete bipartite graph
give a useful Fourier-ratio bound in an explicitly available
labeling? For delivery graphs, can such a bound be derived from
the geometry of the vertex locations and the rule defining an
edge? Can a useful labeling be learned from partial observations
while accounting for the samples used to select it?
\end{question}

It would also be useful to exploit the known bipartition in the
sampling procedure. Can one sample only the $mn$ potentially
present edges, rather than ordered pairs in the full square,
and obtain sharper recovery guarantees? This is particularly
relevant when the part sizes are very different.

For a fixed known labeling with Fourier ratio at most $R$,
how sharp are the dependence on $R^2/\delta$ and the logarithmic
factors in Corollary~\ref{cor:graph-recovery}, which permits at
most $\delta s$ edge errors? Matching lower bounds would clarify
how much can be gained by using symmetry, binary adjacency
values, or additional graph structure. When the labeling must
also be found, the sampling cost should distinguish recovery of
the labeled edge set from recovery up to graph isomorphism.
The obstruction to exact recovery leaves these approximation
questions open.

The bipartite block models of Larremore, Clauset, and Jacobs
\cite{LarremoreClausetJacobs2014} suggest a concrete setting
for these questions. Consider a simple-graph model with
partitions
\[
W=W_1\cup\cdots\cup W_a,
\qquad
S=S_1\cup\cdots\cup S_b.
\]
For each $x\in W_i$ and $y\in S_j$, include the edge
$\{x,y\}$ independently with probability $p_{ij}$. Under a
fixed cyclic labeling, write $f$ for the random adjacency
function and $F=\mathbb E f$ for its mean. Thus $F$ is constant
on each rectangle $W_i\times S_j$ and its transpose, and
vanishes on pairs within either vertex type.

\begin{question}
In parameter ranges where the graph is nonempty with high
probability, which assumptions on the block sizes, probabilities,
and labeled block sets yield useful high-probability bounds
for $\operatorname{FR}(f)$? When do these bounds make
Corollary~\ref{cor:graph-recovery} require substantially fewer
than $mn$ observations?
\end{question}

The distinction between $F$ and $f$ is essential. A small Fourier
ratio for the mean adjacency function does not by itself give
the required bound for a realization. The contribution of the
random fluctuations $f-F$ must also be controlled.

\subsection{Community detection from partial observations}

The known bipartition specifies the two vertex types; it does
not specify the communities within either type. This raises
the question of whether recovering most edges is sufficient
to recover those communities.

Let $B$ and $B_*$ be the $m\times n$ zero-one incidence
matrices of the original and recovered bipartite graphs. Here, rows are indexed by $W$ and columns by $S$, so each
undirected edge corresponds to a single matrix entry.
Unlike the full symmetric adjacency matrix, these rectangular
matrices do not also include the reverse orientation.
Each incorrect edge changes exactly one entry of these
matrices. Consequently, the conclusion of
Corollary~\ref{cor:graph-recovery} gives
\[
\|B-B_*\|_{\mathrm{op}}
\leq \|B-B_*\|_{\mathrm F}
=\sqrt{|E\triangle E_*|}
\leq\sqrt{\delta s},
\]
where $\|\cdot\|_{\mathrm{op}}$ and $\|\cdot\|_{\mathrm F}$
denote the operator and Frobenius norms, respectively.

\begin{question}
For the bipartite block model above, under what assumptions
on community sizes and connection probabilities can a
community-detection method applied to $B_*$ be proved to
recover the underlying communities approximately? More
precisely, can the edge-error bound be converted into a bound
for the number of incorrectly assigned vertices, after
matching community labels within each vertex type? When can
partial observations retain the community-recovery accuracy
obtainable from the full adjacency matrix?
\end{question}

For methods based on singular vectors, the displayed estimate
suggests comparing $\sqrt{\delta s}$ with the relevant
singular-value gap. Control of the resulting subspaces would
then have to be combined with separation conditions that
permit reliable assignment of vertices to communities.
The total edge-error bound allows errors to concentrate on
a small community, so it alone does not ensure accurate
community recovery.

It would also be useful to compare reconstruction followed
by community detection with fitting a block model directly
to the same observed entries. Such a comparison should
account for the observation budget, the accuracy of the
inferred communities, and the computational cost.

\section*{Computational supplement}

The archive \nolinkurl{Fourier_certificates_Oct24_2026.zip}
accompanies this paper. It contains the verifiers for
Propositions~\ref{prop:finite-counterexample}
and~\ref{prop:single-exchange-obstruction}, their shared
arithmetic module, tests, a README with instructions, and
recorded computational results.

The second author developed the original program verifying the
blow-up computation. Because that program was unavailable
for the arXiv submission, the authors prepared
an independent replacement in collaboration with OpenAI's
ChatGPT, following the mathematical specification in
Subsection~\ref{sec:certificate}. The descent verifier,
shared arithmetic module, and tests were prepared by the
authors in the same collaboration. The authors verified the
implementations by examining the source code and running
the programs on their own machines.

The supplied reports and test logs are the original
development-run records; they do not record the authors'
subsequent runs on their own machines.

\appendix
\section{The exact certificate for the four-point obstruction}
\label{app:single-exchange-certificate}

We justify the finite computation in
Proposition~\ref{prop:single-exchange-obstruction}. The method
is the same as in Subsection~\ref{sec:certificate}, with the
parameters chosen for four-element subsets of $\mathbb Z_{17}$.
We give the error estimates explicitly so that the comparison
can be checked independently of the printed decimal values.

Put $Q=10^9$. Starting from Machin's identity
\[
\pi=16\arctan(1/5)-4\arctan(1/239),
\]
truncate each arctangent series after forty terms, indexed from
$0$ to $39$, to obtain a rational number $\pi_0$. The
alternating-series estimate gives
\[
|\pi-\pi_0|\leq e_\pi,
\qquad
e_\pi=\frac{16}{81\cdot5^{81}}
      +\frac{4}{81\cdot239^{81}}.
\]
For $j=0,\ldots,16$, put $\theta_j=2\pi_0j/17$. Let $T_c$
and $T_s$ be the Taylor polynomials of degree at most $79$
for cosine and sine at zero, respectively. Each has forty
nonzero terms; the actual degree of $T_c$ is $78$, since
the coefficient of degree $79$ vanishes. The verifier checks
by exact rational arithmetic that
\[
0\leq\theta_j<7\quad(0\leq j\leq16),
\qquad
2e_\pi+\frac{7^{80}}{80!}<\frac1Q.
\]
Both Taylor remainders are bounded by $7^{80}/80!$.
Moreover, $|\theta_j-2\pi j/17|\leq2e_\pi$, and both sine
and cosine have derivatives bounded in absolute value by one.
Thus each polynomial value approximates the corresponding
trigonometric value at $2\pi j/17$ with error less than $1/Q$.

Define the integer tables
\[
a_j=\lfloor Q T_c(\theta_j)\rfloor,
\qquad
b_j=\lfloor Q T_s(\theta_j)\rfloor.
\]
Taking the floor introduces less than another $1/Q$ of error,
including when the polynomial value is negative. Consequently,
\[
\left|\frac{a_j}{Q}-\cos\frac{2\pi j}{17}\right|<\frac2Q,
\qquad
\left|\frac{b_j}{Q}-\sin\frac{2\pi j}{17}\right|<\frac2Q.
\]

Fix a four-element set $A\subset\mathbb Z_{17}$. For
$1\leq t\leq8$, form
\[
X_t=\sum_{c\in A}a_{tc\bmod17},
\qquad
Y_t=\sum_{c\in A}b_{tc\bmod17},
\qquad
R_t=\left\lfloor\sqrt{X_t^2+Y_t^2}\right\rfloor.
\]
The complex number $(X_t+iY_t)/Q$ approximates
$\overline{S_A(t)}$, whose magnitude equals $|S_A(t)|$,
with error less than $8\sqrt2/Q<12/Q$. The integer square
root is checked by the exact inequalities
\[
R_t^2\leq X_t^2+Y_t^2<(R_t+1)^2.
\]
The triangle inequality and these bounds therefore give
\[
\frac{\max\{0,R_t-12\}}Q
\leq |S_A(t)|
\leq\frac{R_t+13}Q.
\]
The additional $1$ in the upper numerator accounts for the
integer square-root floor.

Since $|S_A(17-t)|=|S_A(t)|$, the nonzero frequencies form
eight conjugate pairs. Hence, for $\Phi(A)=D_{17}(A)/17$,
\[
L(A)\leq\Phi(A)\leq U(A),
\]
where
\[
L(A)=\frac{2}{17Q}\sum_{t=1}^{8}\max\{0,R_t-12\},
\qquad
U(A)=\frac{2}{17Q}\sum_{t=1}^{8}(R_t+13).
\]
Here the integers $R_t$ are calculated separately for each set
$A$. Conjugate symmetry is used for the true Fourier sums;
no exact symmetry of the rounded tables is required.

For $C=\{0,1,4,5\}$ and $P=\{0,1,2,3\}$, the verifier
constructs all $52$ distinct sets in $\mathcal N(C)$ and
checks the exact rational inequalities
\[
U(P)<L(C),
\qquad
U(C)<L(C')\quad\text{for every }C'\in\mathcal N(C).
\]
It also checks the four decimal thresholds in
Proposition~\ref{prop:single-exchange-obstruction}, interpreted
as rational numbers. The minimum over the neighbors is enclosed by
\[
\min_{C'\in\mathcal N(C)}L(C')
\leq\min_{C'\in\mathcal N(C)}\Phi(C')
\leq\min_{C'\in\mathcal N(C)}U(C').
\]
The two endpoint minima need not be attained by the same neighbor.
Only after the exact comparisons have passed are the lower and
upper endpoints rounded downward and upward, respectively, to
produce the displayed intervals.

The accompanying computational supplement supplies the verifier
and its shared arithmetic module:
\begin{center}
\texttt{verify\_descent\_obstruction.py}\\
\texttt{certified\_fourier.py}
\end{center}
All proof-relevant calculations use arbitrary-precision integers
and exact rational arithmetic; no floating-point comparison is
used to certify an inequality.

\end{document}